%% file: Simon,Behrens_multiscale_2018-02-21.tex
\documentclass[10pt,a4paper,footsepline,headsepline,headinclude]{scrartcl}

\usepackage{setspace}

\usepackage{amsmath,amssymb,amsthm,amsfonts,mathrsfs}
\usepackage{a4wide}
\usepackage{scrpage2}

\usepackage[english]{babel}
\usepackage[T1]{fontenc}
\usepackage[utf8x]{inputenc}

\setfootsepline{.5pt}

\usepackage{graphicx}
\usepackage{color}
\usepackage{pdfpages}

\graphicspath{{./Figures/}} 
\usepackage[hangindent=0pt,singlelinecheck=on,format=plain,font=scriptsize]{caption}
\usepackage[font=scriptsize,labelformat=parens]{subcaption}

\usepackage{float}
\restylefloat{figure}
\usepackage{wrapfig}

\theoremstyle{plain}

\theoremstyle{definition}

\theoremstyle{remark}

\newcommand {\peclet}{\mathrm{Pe}}

\newcommand {\eps} {\varepsilon}
\newcommand {\vphi} {\varphi}
\newcommand {\bvphi} {\boldsymbol\varphi}

\DeclareMathOperator{\linhull}{span}

\newcommand{\set}[1]{\left\lbrace #1 \right\rbrace}
\newcommand{\ve}[1]{\mathbf{#1}}

\newcommand{\norm}[2]{\left\| #1 \right\|_{#2}}

\renewcommand{\d}{\:\mathrm{d}}

\lohead{K. Simon, J. Behrens}
\rohead{MsFEM for transient ADEs}
\rofoot{\pagemark} \lofoot{} 

\begin{document}

\begin{center}
  \LARGE \textbf{Multiscale finite elements through advection-induced
    coordinates for transient advection-diffusion equations}
\end{center}

\vspace*{0.2cm}

\begin{center}
  Konrad Simon\textsuperscript{$\ast$} and J\"orn Behrens\textsuperscript{$\ast$}
\end{center}

\vspace*{0.2cm}

\begin{center}
  \small {\textsuperscript{$\ast$}University of Hamburg, Department of Mathematics, \\Grindelberg 5, 20144 Hamburg, Germany \\[0.5cm]
  }
\end{center}

\thispagestyle{empty}

\vspace*{0.5cm}

\begin{abstract}
  \textbf{Abstract.} \small Long simulation times in climate sciences
  typically require coarse grids due to computational
  constraints. Nonetheless, unresolved subscale information
  significantly influences the prognostic variables and can not be
  neglected for reliable long term simulations. This is typically done
  via parametrizations but their coupling to the coarse grid variables
  often involves simple heuristics. We explore a novel up-scaling
  approach inspired by multi-scale finite element methods. These
  methods are well established in porous media applications, where
  mostly stationary or quasi stationary situations prevail. In
  advection-dominated problems arising in climate simulations the
  approach needs to be adjusted. We do so by performing coordinate
  transforms that make the effect of transport milder in the vicinity
  of coarse element boundaries. The idea of our method is quite
  general and we demonstrate it as a proof-of-concept on a
  one-dimensional passive advection-diffusion equation with
  oscillatory background velocity and diffusion.
\end{abstract}

\vspace*{0.5cm}

\input{sec-1.tex}
\input{sec-2.tex}
\input{sec-3.tex}

\input{sec-4.tex}

\paragraph{Acknowledgements.} This work was supported by German
Federal Ministry of Education and Research (BMBF) as Research for
Sustainability initiative (FONA); \textit{www.fona.de} through Palmod
project (FKZ: 01LP1513A).

\bibliography{reference}{}
\bibliographystyle{plain}

\end{document}

%% file: sec-1.tex
\section{Introduction}\label{s-1}

\paragraph{Motivation and Overview.} Geophysical processes in our
atmosphere and in the oceans involve many different spatial and
temporal scales. Reliable climate simulations hence need to take into
account the interaction of many scales since these are coupled and
influence each other. Resolving all relevant scales and their
interactions poses immense computational requirements. Even on modern
high-performance computers computational constraints force us to make
compromises, i.e., to use a limited grid resolution which then often
neglects certain scale interactions or to model interactions in a
heuristic manner.

For example long-term climate simulations as done in the simulation of
paleo climate use grid resolutions of approximately 200km and
more. This of course ignores fine scale processes that can not be
resolved by this resolution although they significantly influence
prognostic variables on the coarse grid. Examples include (but are not
limited to) moving ice-shields, land-sea boundaries, flow over rough
orography, cloud physics and precipitation. None of these processes is
resolved by the grid and current climate simulations cope with this by
using so-called parametrizations. These can be seen as replacements or
simplifications of subgrid processes or processes that are too complex
to be taken into account on the prognostic (coarse) scale. Coupling to
the coarse grid is often done heuristically and can even cause
convergence issues when refining the computational grid.

It is our aim to improve the process of transferring information from
the subgrid scale to the coarse grid in a mathematically consistent
way. Numerical multiscale modelling offers promising mathematical
frameworks to achieve this goal. Such methods are already quite
established in other communities such as the porous media community
but rarely made their way to climate simulations. Most mentionable
here are methods based on homogenization, the heterogeneous multiscale
method (HMM), multiscale methods based on the Eulerian-Lagrangian
localized adjoint method (ELLAM and MsELLAM), variational multiscale
methods and multiscale finite element methods (MsFEM).

Homogenization theory~\cite{Bensoussan2011, Jikov2012, Pavliotis2008}
was originally developed as an analytical tool to study the behavior
of solutions to equations with coefficients that rapidly oscillate on
fast scales $\mathcal{O}(\eps)$ in the limit $\eps\rightarrow 0$. An
equation with effective coefficients is developed that correctly
describes the influence of fast scales as the scale separation
grows. The process of incorporating fine scale information into
coarser scales is usually referred to as upscaling. This does,
however, not resolve fine scales. Also analytical results for
effective equations are very difficult to obtain and available only in
a limited number of cases.

The HMM introduces a framework for a large range of multiscale
problems~\cite{Weinan2003, Weinan2011, Abdulle2012} rather than a
concrete method. In contrast to many other techniques it constitutes a
top-down approach based on blending very different natures of macro
and micro models. One simulates a usually incomplete macro-scale model
by choosing an appropriate method (e.g., a finite element method when
dealing with a variational problem). Missing data is incorporated
wherever it is needed by performing constrained micro-scale
simulations. The micro-scale data then have to be incorporated into a
properly chosen macro solver. The HMM offers a very general approach
but needs a careful choice of the macro and micro models and their
interaction, i.e., the compression and reconstruction of information.

MsELLAM~\cite{Cheng2010, Wang2009} is based on ELLAM~\cite{Celia1990,
  Herrera1993, Russell2002} which uses a space time finite element
method (FEM). Basis functions that satisfy the adjoint equation are
constructed locally in time, which leads to the elimination of
non-boundary terms in the weak form. Back-tracking characteristics
makes it possible to deal with advection dominated problems but leads
to discontinuities in time.

Variational multiscale methods first emerged in~\cite{Brezzi2000,
  Hughes1995, Hughes1998}. The idea here is to decompose the solution
space into a coarse and a fine scale part and to then separate the
usual variational form into a part that is tested with coarse scale
test functions and a part that is tested with fine scale test
functions. This approach, in contrast to the HMM and to
homogenization, resolves fine scale features of the numerical
solution. Mixed variational methods have later been introduced
in~\cite{Arbogast2000, Arbogast2006}, see also~\cite{Graham2012}

The MsFEM is closely related to the variational multiscale method and
has originally been introduced in~\cite{Babuska1983, Babuska1994} and
later been substantially refined in~\cite{Hou1997, Hou1999}. It relies
on the use of non-polynomial basis functions that reflect the fine
scale behavior of the solution. Basis functions are usually
constructed by requiring them to satisfy the equation to leading
order. Many variations have been introduced since then,
see~\cite{Efendiev2009} for a survey. The method is very attractive
since it allows for massive parallelization and is therefore very
scalable.

All the above mentioned methods work well for elliptic problems and
stand in contrast to classical adaptive mesh refinement techniques
(AMR) which essentially constitute a down-scaling approach rather than
an upscaling framework by locally resolving ``interesting'' regions of
the grid~\cite{Jablonowski2004, Behrens2006}. They are proven to work
well also for hyperbolic problems. Our list of methods is by far not
complete. We therefore refer in particular to more comprehensive texts
on different parts of multiscale modeling~\cite{Behrens2006,
  Graham2012, Efendiev2009, Weinan2003, Bensoussan2011} and to
references therein.

\paragraph{Contribution.} In this work we introduce a generalization
of the classical MsFEM introduced in~\cite{Hou1997} to suggest an idea
for dealing with problems that arise in long-term simulations of
climate. The partial differential equations used for these simulations
are mostly time-dependent and are dominated by advection. We therefore
demonstrate our idea in a simple setting on a one-dimensional
advection diffusion equation on a periodic domain that is dominated by
the advective term with fine scale diffusion. In such a scenario the
classical MsFEM which is designed for elliptic problems fails since it
relies on a decomposition of the computational domain into a number of
coarse blocks. These coarse blocks are the support of the modified
basis functions but blocks are essentially decoupled from each
other. Hence information travelling through the entire domain is
blocked at the coarse block boundaries. This means that boundary
conditions of the basis functions need to be chosen with care. Also
the basis functions need to be transient themselves.

Choosing appropriate boundary conditions for the local problems (here
for the basis functions) is a general problem in multiscale
computations. Posing Dirichlet boundary conditions as in the classical
MsFEM on the coarse block boundaries basis functions can develop
boundary layers that do not represent any features that are apparent
in the actual flow. We intend to circumvent this problem not by a
different choice of boundary conditions but by posing Dirichlet
conditions on curves that account for the advective part of the
equation.

Our key idea is based on the fact that for simple background velocity
fields one can transform an advection-diffusion equation into a pure
diffusion equation by following the characteristics of the advective
part, i.e., we move to a Lagrangian setting where we only ``see''
diffusion across the characteristics. Such a transform can be done
without any danger of crossing characteristics if the background
velocity does not depend on space (it may depend on time though). If
the background velocity depends on space we either average in space
and follow the mean flow characteristics or we follow the
characteristics at the coarse grid boundaries and interpolate them
inside the blocks. Each procedure introduces a transform of the
equation that reduces the influence of the advective term.

For a well-behaved background velocity, for example taken from a
coarse grid in a climate simulation, characteristics do not come too
close. Inside each coarse block interpolation of the characteristics
emerging at the coarse grid boundary ensures that (interpolated)
characteristics do not cross or come too close. This amounts to a
setting that is not fully Lagrangian but still makes the advective
term milder. The transformed equation then is very close to a pure
parabolic setting inside a coarse block and purely parabolic in the
vicinity of the boundary of the block. The latter ensures that no
(advective) features cross the boundary and the coarse grid blocks are
not only computationally but also physically decoupled. Posing
Dirichlet boundary conditions on coarse blocks that move with the flow
is therefore useful since information flow is not artificially blocked
as with fixed coarse blocks.

From the computational point of view our method is quite
attractive. It is composed of two parts: an offline phase that
precomputes the basis functions and an online phase that uses these
basis functions to compute the coarse solution. The overhead of
precomputing the basis functions in each coarse block can further be
reduced by parallelization as in the classical elliptic MsFEM. The
online phase is approximately as fast as a low resolution standard FEM
but still reveals fine scale features of a highly resolved solution
and is therefore much more accurate than a standard FEM.

Also note that although we suggest to use the precomputed basis
functions in a finite element framework that the idea is much more
generic. Such modified basis functions can potentially be used in a
different global framework that uses a finite volume or discontinuous
Galerkin formulation.

We are of course aware of the fact that methods for a simple model
like the advection-diffusion equation do not immediately generalize to
more important and harder problems like the primitive
equations. Nevertheless, we strive to suggest ideas for dealing with
the interaction of different parts of such simulations, i.e., between
advective and (parametrized) diffusive features that have a multiscale
character. Besides that our method can potentially be useful for
simpler problems dealing with passive tracer transport. We also point
out problems with this approach and with its generalization. At this
point we are unfortunately not able to give a profound analysis of the
method but we show test cases that will demonstrate the superiority of
our modified MsFEM over standard methods in certain flow regimes. A
rigourous analysis is left for future work.

\paragraph{Outline.} This work is organized as follows. In
Section~\ref{s-2} we introduce the modified MsFEM method. This is done
systematically starting from an easy setting to the general
one. Section~\ref{s-3} shows tests, each of them revealing a different
aspect of the modified MsFEM. We conclude with a discussion of
strengths and limitations of our method in Section~\ref{s-4}.

%% file: sec-2.tex
\section{Description of the Method}\label{s-2}

In this section we will give an overview of our method which is based
on multiscale finite element methods. We first introduce the general
idea of time-dependent basis functions and formulate the discrete
equations. Since this idea does not work in general we then propose
transforms of the original equation and reformulate the discretized
equation. We develop these transforms hierarchically to provide the
reader with a better understanding for the method.

Multiscale finite element methods for elliptic problems have been
intensively investigated by the porous media
community~\cite{Efendiev2009, Graham2012}. In contrast to the
HMM~\cite{Weinan2003} they are designed not only to represent the
coarse scale features of the solution to the problem correctly but
also reveal the fine structure of the solution since the fine scale
behavior is resolved by the basis functions. We strive to adopt this
approach to one-dimensional transient advection-diffusion equations of
the form
\begin{equation}
  \label{eq:2-1}
  \begin{split}
    \partial_t u(x,t) + c(x,t)\partial_x u(x,t) & = \partial_x( \mu(x,t)\partial_x u(x,t)  ) + g(x)\:, \quad (x, t)\in I\times (0,T] \\
    u(0,t) & = u(1, t) \;, \hspace{2.5cm} t\in (0,T] \\
    u(x,0) & = f(x) \:, \hspace{3cm} x\in I
  \end{split}
\end{equation}
on the unit interval $I=[0,1]$ in a periodic setting with rapidly
oscillating coefficients
$c, \mu \in L^\infty([0,T], L_{\text{per}}^\infty(I))$, and smooth
periodic $f,g\in L_{\text{per}}^2(I)$. In order for~(\ref{eq:2-1}) to
be well-posed we assume that $\mu>0$. Then, by the theory of parabolic
equations~\cite{Evans10} the weak form of problem~(\ref{eq:2-1}),
i.e., find $u$ with
\begin{multline}
  \label{eq:2-1a}
  \frac{\d}{\d t}\int_I u(x,t)\vphi(x) \d x + \int_I c(x,t) \partial_x u(x,t) \vphi(x) \d x \\ 
  = -\int_I \mu(x,t) \partial_x u(x,t) \partial_x\vphi(x) \d x + \int_I g(x)\vphi(x) \d x \quad \forall \vphi\in H_{\text{per}}^1(I) \:, t\in (0,T] \:,
\end{multline}
has a unique solution $u$ that is in $L^2([0,T], H_{\text{per}}^1(I))$
and its derivative $\partial_t u$ is in
$L^2([0,T], H_{\text{per}}^{-1}(I))$. By compactness we additionally
have $u\in C([0,T], L_{\text{per}}^2(I))$ and hence the initial
condition $u(x,0) = f(x)$ makes sense. We will impose further
assumptions on the data in the sequel.

\paragraph{Regime of the Data.} Motivated by typical problems that
arise in the parametrization of subgrid processes in climate
simulations we assume that $\mu=\mu(x/\eps, t)$ and
$c=c(x,x/\delta, t)$ with $\eps\ll H$ and $\delta\ll H$ or
$\delta\gtrsim H$, where $H$ denotes a the magnitude of a resolved
numerical scale. By this we mean that oscillations in the data can not
be resolved by $H$, i.e., one can even think of $\epsilon,\delta$ as
vectors that represent regimes of scales in the data. However, we
suppress the $\epsilon, \delta$ notation unless necessary
otherwise. The reader should note that although we keep the idea of
possibly fine oscillatory scales $O(\delta)$ in the background
velocity the standard case is that $\delta\gtrsim H$. This is because
background velocity data usually comes from the coarse simulation
scale whereas $\mu$ comes from parametrizations and therefore lives on
the subgrid scale. However, to be as general as possible we would like
to keep both options, i.e., $\delta$ to denote resolved or unresolved
scales.

Furthermore, we are interested in a regime that is dominated by the
advective term. Since the coefficients are assumed to be oscillatory
it is less enlightening to use a single P\'eclet number
$\peclet = \norm{c}{L^\infty}L/\mu$, where $L$ is a characteristic
length scale, since in our case $c$ can be is a very local quantity as
well. Instead we are interested in the regime where
$\peclet(x) = c(x,t)L/\mu(x,t)$ is large in average, i.e., we suggest
to look at $\peclet$ with a large (normalized) $L^1$-norm. This would
mean that the tracer distribution over time is on average dominated by
advection, i.e., for the background velocity we assume that it is
dominated by its mean.

\paragraph{The Basic Idea of the Method.} Multiscale finite elements
consist of two main components. A global formulation and modified
basis functions, i.e., a localization. In our approach we decompose
the domain $I$ into a number coarse cells $K_i$, $i=1,\dots,N-1$ of
mesh width $H$. With each Eulerian node $x_m$, $m=1,\dots,N$, of the
coarse mesh we associate a basis function $\vphi_m^{\mathrm{ms}}$ that
satisfies $\vphi_m^{\mathrm{ms}}(x_n)=\delta_{mn}$, i.e., globally we
use a conformal finite element formulation. Standard finite elements
interpolate the basis functions between the nodes in a prescribed
functional way, often polynomial. For the multiscale finite element
basis functions used in our global formulation, i.e., the second
constituent, we interpolate in such a way that the basis contains
information about the fine scale structure of the problem at hand. In
the porous media community the common way is to require the basis to
satisfy the homogeneous equation locally, i.e., in each coarse grid
block, to the highest order. We aim to adopt this idea to transient
equations that contain advective terms which is the major difficulty.

In each coarse block $K_i$ we then seek basis functions $\vphi_i^l$,
$l=1,2$, that satisfy
\begin{equation}
  \label{eq:2-2}
  \begin{split}
    \partial_t \vphi_i^l(x,t) + c(x,t)\partial_x \vphi_i^l(x,t) & = \partial_x( \mu(x,t)\partial_x \vphi_i^l(x,t)  ) \:, \quad (x, t)\in K_i\times (0,T] \\
    \vphi_i^l(0,t) & = \vphi_{i,0}^l(0) \;, \hspace{2.5cm} t\in (0,T] \\
    \vphi_i^l(1,t) & = \vphi_{i,0}^l(1) \;, \hspace{2.5cm} t\in (0,T] \\
    \vphi_i^l(x,0) & = \vphi_{i,0}^l(x) \:, \hspace{3cm} x\in K_i
  \end{split}
\end{equation}
in the weak sense where $\vphi_{i,0}^l$ is the $l$-th $P_1$-finite
element basis function on the coarse scale. Note that these problems
have unique weak solutions
$\vphi_i^l\in L^2([0,T],\vphi_{i,0}^l+H_0^1(I))$ with
$\partial_t\vphi_i^l\in L^2([0,T], H^{-1}(I))$ and hence
$\vphi_i^l\in C([0,T],L^2(I))$. To numerically solve these problems we
employ the method of lines, i.e., we discretize in space and then
solve the resulting ODE. Using a conformal standard FEM in space the
problem reads: Find $\vphi_i^l\in C^1([0,T], V^h(K_i))$ such that
\begin{multline}
  \label{eq:2-2a}
    \int_I \partial_t \vphi_i^{l,h}(x,t)\psi(x) \d x + \int_I c(x,t) \partial_x \vphi_i^{l,h}(x,t) \psi(x) \d x \\ 
  = -\int_I \mu(x,t) \partial_x\vphi_i^{l,h}(x,t) \partial_x \psi(x) \d x \quad \forall \psi\in V^h(K_i) \:, t\in (0,T] 
\end{multline}
with the initial condition
$\phi_i^{l,h}(x,0) = P_{V^h(K_i)}\vphi_{i,0}^l(x)$ where $P_{V^h}$ is
the projection onto $V^h(K_i)$. The space $V^h(K_i)$ is given by
\begin{equation}
  \label{eq:2-2b}
  V^h(K_i) = \set{ \psi \in H^1(K_i) \: | \: \text{$\psi$ is piece-wise polynomial and continuous}  } \:.
\end{equation}

The nodal basis functions $\vphi_m$ are then constructed by ``gluing''
the complementing parts of two basis functions of two adjacent cells,
i.e., the basis functions have the same boundary conditions at the
intersection node like in a standard FEM. For an illustration see
Figure~\ref{fig-basis}. \input{fig-basis}

The basis functions need to be pre-computed in an offline-phase. This
can be done in parallel for each cell since the local problems do not
depend on each other. In each of the small local problems actually
only one basis function $\vphi_i^l$ has to be computed since the other
(unique) solution can be computed as
$\vphi_i^k(x,t) = 1 - \vphi_i^l(x,t)$, $k\not= l$, due to
linearity. Consequently the basis functions fulfill
$\sum_{j=1}^{N} \vphi_j(x,t) = 1$, i.e., they constitute a
decomposition of unity.

For the global formulation, in contrast to standard finite elements
where basis functions do not depend on time, we seek solutions of the
form
\begin{equation}
  \label{eq:2-3}
  u^H(x,t) = u_j^H(t)\vphi_j(x,t) \:.
\end{equation}
Here we used the Einstein's sum convention. With this ansatz we design
our global weak form as: Find $u^H\in V_{\text{per}}^H$ such
that
\begin{multline}
  \label{eq:2-4a}
  \int_I \partial_t u^H(x,t)\vphi(x,t) \d x + \int_I c(x,t) \partial_x u^H(x,t) \vphi(x,t) \d x \\ 
  = -\int_I \mu(x,t) \partial_x u^H(x,t) \partial_x\vphi(x,t) \d x \quad \forall \vphi\in V_{\text{per}}^H(t) \:, t\in (0,T] 
\end{multline}
with the initial condition $u^H(x,0) = P_{V_{\text{per}}^H (0)}f$
where $P_{V_{\text{per}}^H (0)}$ denotes the projection onto the
subspace $V_{\text{per}}^H(0)\subset L_{\text{per}}^2(I)$. Using the
method of lines the global finite element space $V_{\text{per}}^H$ can
be formulated as
\begin{equation}
  \label{eq:2-4b}
  V_{\text{per}}^H = \set{u^H\in C^1([0,T], H_{\text{per}}^1(I)) \: | \: u^H(t) \in V_{\text{per}}^H(t)}
\end{equation}
with
\begin{equation}
  \label{eq:2-4c}
  V_{\text{per}}^H(t) = \linhull_m\set{\vphi_m(\cdot ,t) \: | \: \vphi_m \: \text{satisfies the discretized equation~(\ref{eq:2-2}) at time}\: t \in[0,T]}
\end{equation}
Note that a basis function $\vphi_m$ is glued together from two parts
that are element of the space $V^h(K_i)$. This renders the method
conformal in space. We use the same constructed basis as basis of the
test function space. However, a Petrov-Galerkin formulation using,
e.g., standard $P_1$-basis functions is possible.

Using expression~(\ref{eq:2-3}) as ansatz for the solution
in~(\ref{eq:2-4a}) we get
\begin{eqnarray}
  \label{eq:2-4d}
  \int_I \vphi_i(x,t) \partial_t \left( u_j^H(t) \vphi_j(x,t)\right) \d x + \int_I \vphi_i(x,t) c(x,t) u_j^H(t)\partial_x\vphi_j(x,t) \d x \\ 
  = -\int_I \partial_x\vphi_i(x,t) \mu(x,t) u_j^H(t)\partial_x\vphi_j(x,t) \d x +\int_I \vphi_i(x,t) g(x) \d x \quad \forall i=1,\dots,N \:, t\in (0,T] \:.
\end{eqnarray}
which can be expressed as an ODE:
\begin{equation}
  \label{eq:2-5}
  M(t)\frac{\d}{\d t}\ve{u}^H(t) + N(t)\ve{u}^H(t) + A(t)\ve{u}^H(t) = D(t)\ve{u}^H(t) + \ve G(t) \:,
\end{equation}
where $\ve{u}^H(t) = (u_1^H(t), \dots, u_N^H(t))^T$,
$\ve G(t)_i = \int_I \vphi_i(x,t) g(x) \d x$ and
\begin{equation}
  \label{eq:2-6}
  \begin{split}
    M(t)_{ij} = \int_I \vphi_i(x,t)\vphi_j(x,t) \d x \:, & \quad N(t)_{ij} = \int_I \vphi_i(x,t) \partial_t\vphi_j(x,t) \d x \:, \\
    A(t)_{ij} = \int_I \vphi_i(x,t) c(x,t) \partial_x\vphi_j(x,t) \d x \:, &
    \quad\text{and} \quad D(t)_{ij} = -\int_I \partial_x\vphi_i(x,t)
    \mu(x,t) \partial_x\vphi_j(x,t) \d x 
  \end{split}
\end{equation}
and with the initial condition
$\ve{u}(0)_i=\int_I\vphi(x,0)f(x) \d x$. Note that the second term
$N(t)$ in~(\ref{eq:2-5}) appears since the basis depends on time. Also
note that this is a purely Eulerian setting.

Now, one could simply solve this ODE using a suitable time
integrator. Unfortunately, this method fails due to the dominance of
the advective term. Simply taking small time steps does not remedy the
problem since the problem lies in the way we compute the solutions to
the local problems. The basis functions are supposed to carry
information about the fine scale structure of the solution but note
that the local problems are posed with Dirichlet boundary
conditions. Such type of boundary condition essentially blocks
advected information at coarse cell boundaries since all coarse cells
are decoupled from each other. Hence, for flow dominated by advection
multiscale basis functions constructed according to~(\ref{eq:2-2})
would develop a boundary layer that is not present in the global
flow. The "correct" choice of type of boundary conditions is indeed a
general problem that one is faced with when solving such local
problems~\cite{Efendiev2009}.

\paragraph{A Modified Idea.} We try to circumvent the problem above by
moving to a different set of coordinates that make the advective term
milder, e.g., we could follow a suitable set of curves emerging at our
coarse grid points as done in a full Lagrangian framework. Then
Dirichlet boundary conditions could be posed on these curves since
these naturally decouple flow regions. Intuitively this would give
better results since the transform will bring us ``closer'' to a
parabolic setting. However, a full Lagrangian framework is, in
general, difficult to handle since it is not easy to implement and
since characteristics could potentially come very close (or even
cross) which requires a special treatment in such regions. We
therefore propose a simpler method.

We start by explaining our idea in a simple setting and then
generalize. For this let the background velocity $c(x,t)=c(t)$, i.e.,
we consider a flow without shear. This special assumption actually
allows us to move to a full Lagrangian setting since we can simply
follow characteristics and know they will not intersect. The new
coordinates $(\xi,\tau)$ are then (implicitly) given by
\begin{equation}
  \label{eq:2-7a}
  \begin{split}
    x(\xi,\tau) & = \xi + \int_0^\tau c(s) \d s \\
    t(\xi,\tau) & = \tau \:.
  \end{split}
\end{equation}
A little bit more general is the choice
\begin{equation}
  \label{eq:2-7b}
  \begin{split}
    x(\xi,\tau) & = \xi + \int_0^\tau \langle c \rangle (s) \d s \:, \quad \text{where} \quad \langle c \rangle (t) = \frac{1}{|I|}\int_I c(x,t) \d x \\
    t(\xi,\tau) & = \tau
  \end{split}
\end{equation}
if $c$ explicitly depends on space. In the new
coordinates~(\ref{eq:2-7b}) equation~(\ref{eq:2-1}) reduces to
\begin{equation}
  \label{eq:2-8}
  \begin{split}
    \partial_\tau \hat u(\xi,\tau) + \tilde c(\xi,\tau) \partial_\xi \hat u(\xi,\tau)
    & = \partial_\xi( \hat\mu(\xi,\tau)\partial_\xi \hat u(\xi,\tau)  ) + \hat g(\xi)\:, \quad (\xi, \tau)\in I\times (0,T] \\
    \hat u(0,\tau) & = \hat u(1, \tau) \;, \hspace{2.5cm} \tau\in (0,T] \\
    \hat u(\xi,0) & = \hat f(\xi) \:, \hspace{3cm} \xi\in I 
  \end{split}
\end{equation}
where $\hat u(\xi,\tau) = u(x(\xi,\tau),t(\tau))$, respectively
$\hat c, \hat \mu, \hat f, \hat g$, and
\begin{equation}
  \label{eq:2-9}
  \tilde c(\xi,\tau) = \hat c(\xi,\tau) - \langle c \rangle (\tau) \:.
\end{equation}
Note that if $c$ is constant or only depends on time, i.e., we
employ~(\ref{eq:2-7a}), then $\tilde c(\xi,\tau)=0$ and~(\ref{eq:2-8})
reduces to a simple diffusion equation. An illustration of the two
transforms is done in Figure~\ref{fig-characteristics}(a)~and~(b).

\input{fig-characteristics} The subscale problem for the basis
functions~(\ref{eq:2-2}) needs to be transformed as well: Find a
(weak) solution $\hat\vphi_i^l(\xi,\tau)$, $l=1,2$, such that
\begin{equation}
  \label{eq:2-10}
  \begin{split}
    \partial_\tau \hat \vphi_i^l(\xi,\tau) + \tilde c(\xi,\tau)\partial_\xi \hat \vphi_i^l(\xi,\tau) & = 
    \partial_\xi( \hat \mu(\xi,\tau)\partial_\xi \hat \vphi_i^l(\xi,\tau)  ) \:, \quad (\xi,\tau)\in \hat K_i\times (0,T] \\
    \hat \vphi_i^l(0,\tau) & = \vphi_{i,0}^l(0) \;, \hspace{2.5cm} \tau\in (0,T] \\
    \hat \vphi_i^l(1,\tau) & = \vphi_{i,0}^l(1) \;, \hspace{2.5cm} \tau\in (0,T] \\
    \hat \vphi_i^l(\xi,0) & = \vphi_{i,0}^l(\xi) \:, \hspace{3cm}
    \xi\in \hat K_i \:.
  \end{split}
\end{equation}
A cell $\hat K_i$ then essentially moves with the flow described by
the transform without being deformed.

\paragraph{A Generalized Modified Idea.} In case of a large scale
separation between the mean flow and the oscillatory parts of the
background velocity the method using transform~(\ref{eq:2-7b})
delivers good results. However, as mentioned earlier, in many
practically relevant scenarios in climate simulations the background
velocity does not exhibit such a scale separation, i.e., there are
mostly (globally resolved) scales $\delta\gtrsim H$.  We would
nevertheless like to keep possible small scale variations in $c$ to be
as general as possible. Variations at a scale $O(H)$ and larger would
have the consequence that the average velocity $\tilde c$ in each cell
$\hat K_i$ can potentially be large in contrast to the case where one
has a scale separation. Consequently basis functions following the
evolution given by~(\ref{eq:2-10}) can still exhibit steep boundary
layers that are not apparent in the actual flow. As a remedy one can
make an effort to stay ``close'' to the full Lagrangian setting
locally, i.e., we can track the characteristics (particle
trajectories) of the flow starting at the coarse nodes. Boundary
conditions for the multiscale basis functions can then be posed on
these characteristics. In scenarios in which characteristics do not
come too close or diverge too strongly this is, as we will show in
numerical examples, a reasonable strategy. Fortunately, in many
situations in climate simulations there is no situation in which true
shocks occur, e.g., wind field data that is taken from a coarse grid
does not behave too badly.

To formalize this strategy let $x(\xi,\tau)$ denote the
characteristics emerging from $\xi \in I$. The evolution is governed
by the set of ODEs that reads
\begin{equation}
  \label{eq:2-11}
  \begin{split}
    \frac{\d}{\d \tau} x(\xi,\tau) & = c(x(\xi,\tau),\tau) \:, \quad x(\xi,0) = \xi \\
    \frac{\d}{\d \tau} t(\xi,\tau) & = 1  \:, \quad t(\xi,0) = 0 \:.
  \end{split}
\end{equation}
This evolution induces a differentiable transform between
$(x,t)$-coordinates and $(\xi,\tau)$-coordinates as long as
characteristics do not intersect. The latter is the case if $c(x,t)$
is locally Lipschitz-continuous in $x$ and continuous in $t$. If we
were to use transform~(\ref{eq:2-11}) everywhere in our domain (what
we do not pursue) this would simply mean a change from Eulerian to
Lagrangian variables. The following computations are formal and only
hold as long as trajectories do not cross. As a first step we will
rewrite the equations in characteristic coordinates. For the gradients
we have the relation
\begin{equation}
  \label{eq:2-12a}
  \nabla_{(x,t)}(\cdot) = \nabla_{(\xi,\tau)} (\cdot)
  \begin{pmatrix}
    \frac{\partial \xi}{\partial x} & \frac{\partial \xi}{\partial t} \\
    \frac{\partial \tau}{\partial x} & \frac{\partial \tau}{\partial t}
  \end{pmatrix}
\end{equation}
and 
\begin{equation}
  \label{eq:2-12b}
  \nabla_{(\xi,\tau)} (\cdot) = \nabla_{(x,t)}(\cdot)
  \begin{pmatrix}
    \frac{\partial x}{\partial \xi} & \frac{\partial x}{\partial \tau} \\
    \frac{\partial t}{\partial \xi} & \frac{\partial t}{\partial \tau}
  \end{pmatrix}
\end{equation}
Hence, we get
\begin{equation}
  \label{eq:2-13}
  \begin{split}
    \partial_x u(x,t) & = \partial_{\xi}\hat u(\xi,\tau)\frac{\partial \xi}{\partial x} + \partial_{\tau}\hat u(\xi,\tau)\frac{\partial \tau}{\partial x} \\
    \partial_t u(x,t) & = \partial_{\xi}\hat u(\xi,\tau)\frac{\partial \xi}{\partial t} + \partial_{\tau}\hat u(\xi,\tau)\frac{\partial \tau}{\partial t} \:.
  \end{split}
\end{equation}
By means of the inverse function theorem the Jacobian
in~(\ref{eq:2-12b}) is (locally) the inverse of the Jacobian
in~(\ref{eq:2-12a}). Using Kramer's rule the inverse is given by
\begin{equation}
  \label{eq:2-14}
    \begin{pmatrix}
    \frac{\partial \xi}{\partial x} & \frac{\partial \xi}{\partial t} \\
    \frac{\partial \tau}{\partial x} & \frac{\partial \tau}{\partial t}
  \end{pmatrix}
 = \frac{1}{D}
  \begin{pmatrix}
     \frac{\partial t}{\partial \tau} & -\frac{\partial x}{\partial \tau} \\
    -\frac{\partial t}{\partial \xi} & \frac{\partial x}{\partial \xi}
  \end{pmatrix}
\end{equation}
where
\begin{equation}
  \label{eq:2-15}
  D = \frac{\partial x}{\partial\xi}\frac{\partial t}{\partial \tau} - \frac{\partial x}{\partial\tau}\frac{\partial t}{\partial\xi}
\end{equation}
is the determinant of the Jacobian in~(\ref{eq:2-12b}).

Now, using~(\ref{eq:2-11}) in~(\ref{eq:2-13}) we see that 
\begin{equation}
  \label{eq:2-16}
  \partial_t u(x,t) + c(x,t)\partial_x u(x,t) 
  = \partial_\tau\hat u(\xi,\tau) + \left[ \hat c(\xi,\tau) - \frac{\partial x}{\partial \tau}\right] \left(\frac{\partial x}{\partial \xi}\right)^{-1}\hat u(\xi, \tau) \:.
\end{equation}
Looking again at~(\ref{eq:2-11}) we see that
\begin{equation}
  \label{eq:2-17}
  \frac{\partial x}{\partial \tau} = \hat c (\xi,\tau)
\end{equation}
and hence the velocity term in the transformed equation would
vanish. The advantage of such a full Lagrangian setting would be that
one has then a purely parabolic equation as in the case when the
background velocity is constant or solely depends on time. But
transform~(\ref{eq:2-11}) breaks down in the presence of
shocks. However, true shocks can not occur since the overall equation
we are dealing with is parabolic (although dominated by advection) and
strictly speaking there are no characteristics. We would nonetheless
still like to reduce the advective effects and retain the ``nice''
part of the equation. This is where our framework starts to differ
from the full Lagrangian setting.

The idea is now, as mentioned earlier, to pose Dirichlet boundary
conditions on the characteristics in order to solve subscale problems
similar to~(\ref{eq:2-10}). In order to do this we only need the
coarse cell boundaries to be characteristics of the advective
term. Inside each coarse cell we just need to reduce the advective
effects of the original equation. This will then prevent flow across
coarse cell boundaries since near the boundary the transformed
equation will be nearly parabolic. Further, the solution to the
equation inside the coarse cell will reflect all the subscale features
of the solution of the original equation.

The diffusive term transforms to
\begin{equation}
  \label{eq:2-18}
  \partial_x(\mu(x,t)\partial_x u(x,t)) = 
  \partial_\xi\left[
    \hat\mu(\xi,\tau)\partial_\xi \hat u(\xi,\tau)\left(\frac{\partial x}{\partial \xi}\right)^{-1}
  \right] 
  \left(\frac{\partial x}{\partial \xi}\right)^{-1} \:.
\end{equation}

Collecting~(\ref{eq:2-16}) and~(\ref{eq:2-18}) together with the
boundary conditions we are now ready to introduce the transformed
problem for the multiscale basis functions as: By the method of lines
find a (variational) solution
$\hat\vphi_i^l \in C^1([0,T], V^h(\hat K_i))$, $l=1,2$, such that
\begin{equation}
  \label{eq:2-19}
  \begin{split}
    \partial_\tau \hat \vphi_i^l(\xi,\tau) + \tilde
    c(\xi,\tau)\partial_\xi \hat \vphi_i^l(\xi,\tau) & =
    \partial_\xi\left[ \hat \mu(\xi,\tau) \left(\frac{\partial x}{\partial \xi}\right)^{-1} \partial_\xi \hat \vphi_i^l(\xi,\tau)  \right] \left(\frac{\partial x}{\partial \xi}\right)^{-1}\\
    & \qquad + \hat g(\xi,\tau) \qquad, (\xi,\tau)\in \hat K_i\times (0,T] \\
    \hat \vphi_i^l(0,\tau) & = \vphi_{i,0}^l(0) \;, \hspace{2.5cm} \tau\in (0,T] \\
    \hat \vphi_i^l(1,\tau) & = \vphi_{i,0}^l(1) \;, \hspace{2.5cm} \tau\in (0,T] \\
    \hat \vphi_i^l(\xi,0) & = \vphi_{i,0}^l(\xi) \:, \hspace{3cm}
    \xi\in \hat K_i
  \end{split}
\end{equation}
where
\begin{equation}
  \label{eq:2-20}
  \tilde c(\xi,\tau) = \left[\hat c(\xi,\tau) - \frac{\partial x}{\partial \tau}\right]\left(\frac{\partial x}{\partial \xi}\right)^{-1}
\end{equation}
and
\begin{equation}
  \label{eq:2-20a}
  V^h(\hat K_i) = \set{ \psi \in H^1(\hat K_i) \: | \: \text{$\hat \vphi$ is piece-wise polynomial and continuous}  } \:.
\end{equation}

In each of the $N-1$ coarse cells $\hat K_i=[\xi_i, \xi_{i+1}]$ we now
have to compute two additional quantities, i.e.,
$\partial x/\partial \xi$ and $\partial x/\partial \tau$. Computation
of the former is based on linear interpolation in space between the
left and right boundary nodes of $\hat K_i$ which follow the
characteristics starting at $\xi_i$ and $\xi_{i+1}$.  (more details in
the paragraph about the implementation). For the latter, i.e.,
$\partial x/\partial \tau$, we use the fact that
$\partial x/\partial \tau (\xi_i,\tau) = \hat c(\xi_i,\tau)$ at each
coarse node. Inside the cells $\hat K_i$ we can again interpolate
linearly the values $\hat c(\xi_i,\tau)$ and $\hat c(\xi_{i+1},\tau)$.
Note that now $\tilde c(\xi_i,\tau)=0$ at all coarse nodes $\xi_i$
globally in time. This means that characteristics of the transformed
system starting at the coarse nodes are separatrices of the dynamical
system $\mathrm{d}/\mathrm{d}s \xi (s) =\tilde c(\xi(s),s)$ and hence
no flow across coarse cell boundaries is possible in these
coordinates, i.e., the flow is separated by the cells $\hat K_i$.

The global weak form can be formulated as: Find
$\hat u(\xi,\tau) \in \hat V_{\text{per}}^H$ such that
\begin{multline}
  \label{eq:2-21a}
  \int_I \partial_\tau \hat u^H(\xi,\tau)\hat \vphi^H(\xi,\tau) \d \xi + \int_I \tilde c(\xi,\tau) \partial_\xi u^H(\xi,\tau) \hat \vphi^H(\xi,\tau) \d \xi \\ 
  = -\int_I 
\left[
    \hat\mu(\xi,\tau)\partial_\xi \hat u^H(\xi,\tau)\left(\frac{\partial x}{\partial \xi}\right)^{-1}
  \right] 
\partial_\xi\left[
\hat\vphi^H(\xi,\tau)\left(\frac{\partial x}{\partial \xi}\right)^{-1}
\right]
\d \xi \\ + \int_I\hat g(\xi,\tau)\hat\vphi^H(\xi,\tau) \d \xi\quad \forall \hat \vphi^H\in V_{\text{per}}^H(\tau) \:, \tau\in (0,T] 
\end{multline}
where $\tilde c$ is given by~(\ref{eq:2-20}). The initial condition is
$\hat u^H(\xi,0) = P_{\hat V_{\text{per}}^H(0)}f$ where
$P_{\hat V_{\text{per}}^H(0)}$ denotes the projection onto the subspace
$\hat V_{\text{per}}^H(0)\subset L_{\text{per}}^2(I)$. The global
finite element space $\hat V_{\text{per}}^H$ can be formulated as
\begin{equation}
  \label{eq:2-21b}
  \hat V_{\text{per}}^H = \set{\hat u\in C^1([0,T], H_{\text{per}}^1(I)) \: | \: \hat u(\tau) \in \hat V_{\text{per}}^H(\tau)}
\end{equation}
with
\begin{equation}
  \label{eq:2-21c}
  \hat V_{\text{per}}^H(\tau) = \linhull_m\set{\hat \vphi_m \: | \: \hat \vphi_m \: \text{satisfies the discretized equation~(\ref{eq:2-19}) at time}\: \tau \in[0,T]} \:.
\end{equation}
Note that this formulation includes the ones using the
transforms~(\ref{eq:2-7a}) or~(\ref{eq:2-7b}) and that it is not fully
Lagrangian if one uses~(\ref{eq:2-11}).

\paragraph{The Algorithm and some Notes on the Implementation.}

In this paragraph we give a brief overview of steps to be taken for
the implementation of our above described method.

\textit{Step 1.} The first step is to set up a uniform mesh with nodes
$\xi_m$, $m=1,\dots,N$, and cells $\hat K_i$, $i=1,\dots,N-1$ with
mesh size $\epsilon, \delta\gtrsim h$. Now one decides for the
transform. For each node $\xi_m$ either uses ~(\ref{eq:2-7a})
or~(\ref{eq:2-7b}) or one solves the ODE~(\ref{eq:2-11}) on the time
interval $[0,T]$ using a suitable solver. For our examples we used an
adaptive Runge-Kutta-4/5 solver.

\textit{Step 2.} One now needs to solve for the basis in an offline
phase. This can be done in parallel since all local problems are
independent of each other. As indicated earlier we solve
equation~(\ref{eq:2-19}) using a conformal standard $P_1$-finite
element formulation. For this we create in each cell $\hat K_i$ a
uniform mesh with cells $\hat k_j^i$, $i=1,\dots,N_f-1$ and nodes
$\xi_n^m$, $n=1,\dots,N_f$, with $\xi_1^m=\xi_m$ and
$\xi_{N_f}^m=\xi_{m+1}$. In cell $\hat K_i$ the discretized ODE for
the vector $\hat\bvphi _i^{l,H}(\tau)$ of degrees of freedom of the
basis function $\hat \vphi_i^{l,H}(\tau)$ reads
\begin{equation}
  \label{eq:2-22a}
  \hat M^{\hat K_i}\frac{\d}{\d \tau} \hat\bvphi _i^{l,H}(\tau) + \hat A(\tau)^{\hat K_i} \hat\bvphi _i^{l,H}(\tau)
  = \hat D(\tau)^{\hat K_i}\hat\bvphi _i^{l,H}(\tau) \:, \; \hat\bvphi_i^{l,H}(0)=\bvphi_{i,0}^{l,H}
\end{equation}
where $\bvphi_{i,0}^{l,h}$ denotes the nodal values of the standard
basis. The system matrices are given by
\begin{equation}
  \label{eq:2-22b}
  \begin{split}
    \hat M_{kl}^{\hat K_i} & = \int_{\hat K_i} \psi_k^h(\xi)\psi_l^h(\xi) \d \xi \:, \\    
    \hat A(\tau)_{kl}^{\hat K_i} & = \int_{\hat K_i} \psi_k^h(\xi) \tilde c(\xi,\tau) \partial_\xi\psi_l^h(\xi) \d \xi \:, \\
    \hat D(\tau)_{kl}^{\hat K_i} &= -\int_{\hat K_i} \left[\hat \mu(\xi,\tau) \partial_\xi\psi_k^h(\xi) \left(\frac{\partial x}{\partial \xi}\right)^{-1} \right]
    \left[\psi_l^h(\xi)\left(\frac{\partial x}{\partial \xi}\right)^{-1} \right] \d \xi
  \end{split}  
\end{equation}
where $\psi_k^h$ denote the basis of standard that consists of the
common hat functions. Note again that the computation of $A(\tau)$ and
$D(\tau)$ involves terms of the form
\begin{equation}
  \label{eq:2-23}
  \frac{\partial x}{\partial \tau} \quad\text{and}\quad \left(\frac{\partial x}{\partial \xi}\right)^{-1} \:.
\end{equation}

Time stepping is done using a second-order explicit scheme for the
advective term and an (implicit) Crank-Nicolson scheme for the
diffusion. Note that the advection matrix and the diffusion matrix
need to be assembled at each time step. All system matrices
in~(\ref{eq:2-22b}) are stored for each time step (for a reason see
\textrm{Step 3}).

In each cell actually only one basis function $\hat \vphi_i^l$ needs
to be computed. By linearity the other (unique) solution is then given
by $\vphi_i^k(x,t) = 1 - \vphi_i^l(x,t)$, $k\not= l$. Hence the
constructed basis functions globally constitute a decomposition of
unity.

\textit{Step 3.} The global system~(\ref{eq:2-21a}) is solved in an
online phase. The discretized equations form a system of ODEs of size
$(N-1)$-by-$(N-1)$ (note the periodic boundary conditions):
\begin{equation}
  \label{eq:2-24}
  \hat M(\tau)\frac{\d}{\d \tau}\ve{\hat u}^H(\tau) + \hat N(\tau)\ve{\hat u}^H(\tau) + \hat A(\tau)\ve{\hat u}^H(\tau) = \hat D(\tau)\ve{\hat u}^H(\tau) + \hat G(\tau)
\end{equation}
where
$\hat G(\tau) = \int_I \hat\vphi_i^H(\xi,\tau)\hat g(\xi,\tau) \d
\xi$,
\begin{equation}
  \label{eq:2-25}
  \begin{split}
    \hat M(\tau)_{ij} = \int_I \hat\vphi_i^H(\xi,\tau)\hat\vphi_j^H(\xi,\tau) \d \xi \:, & \quad \hat N(\tau)_{ij} = \int_I \hat\vphi_i^H(\xi,\tau) \partial_t\hat\vphi_j^H(\xi,\tau) \d \xi \:, \\
    \hat A(\tau)_{ij} = \int_I \hat\vphi_i^H(\xi,\tau) \tilde c(\xi,\tau) \partial_\xi\hat\vphi_j^H(\xi,\tau) \d \xi \:, &
    \quad\text{and} \quad \hat D(\tau)_{ij} = -\int_I \partial_\xi\hat\vphi_i^H(\xi,\tau)
    \mu(\xi,\tau) \partial_\xi\hat\vphi_j^H(\xi,\tau) \d \xi \:.
  \end{split}
\end{equation}
For the time stepping we again use the same second-order explicit
scheme for the advective term and the (implicit) Crank-Nicolson scheme
for the diffusion that we already used for the basis. Since all the
matrices in~(\ref{eq:2-25}) need to be assembled at each time step
solving the global system can potentially become very expensive. By
using a trick looping over all coarse cells and then over all fine
cells can be circumvented since we stored all system matrices
in~(\ref{eq:2-22b}). To see this note that each global basis function
is a linear combination of fine scale (standard) basis functions:
\begin{equation}
  \label{eq:2-26}
  \hat\vphi_i^H(\xi,\tau) = \alpha_{ij}(\tau)\hat \psi_j^h(\xi)
\end{equation}
As a generic example we show this for the advective term. By virtue
of~(\ref{eq:2-26}) we have for the global element (advection) matrix
for the element $\hat K_n$
\begin{equation}
  \label{eq:2-27}
  \begin{split}
    \hat A(\tau)_{kl}^n & = \left( \int_{\hat K_n} \vphi_i^H(\xi,\tau) \tilde c(\xi,\tau) \partial_\xi\vphi_j^H(\xi,\tau) \d \xi \right)_{kl} \\
    & = \int_{\hat K_n} \alpha_{ik}(\tau)\hat \psi_k^h(\xi) \tilde c(\xi,\tau) \partial_\xi \alpha_{jl}(\tau)\hat \psi_l^h(\xi) \d \xi \\
    & = \alpha \hat A^{\hat K_n}(\tau) \alpha^T
  \end{split}
\end{equation}
where $\hat A^{\hat K_n}(\tau)$ is the stored (advection) matrix of
the coarse cell $\hat K_n$ given by~(\ref{eq:2-22b}). Hence,
assembling the global system at each time step breaks down in
collecting information from each coarse block and a simple matrix
multiplication. This renders the online phase very fast. Also, to
ensure the unique solvability of the characteristic equation we
additionally need to assume that $c(x,\cdot)$ be continuous and that
$c(\cdot,t)$ be Lipschitz-continuous.

%% file: fig-basis.tex
\begin{figure}[b!]
  \centering
  \includegraphics[width=0.95\textwidth]{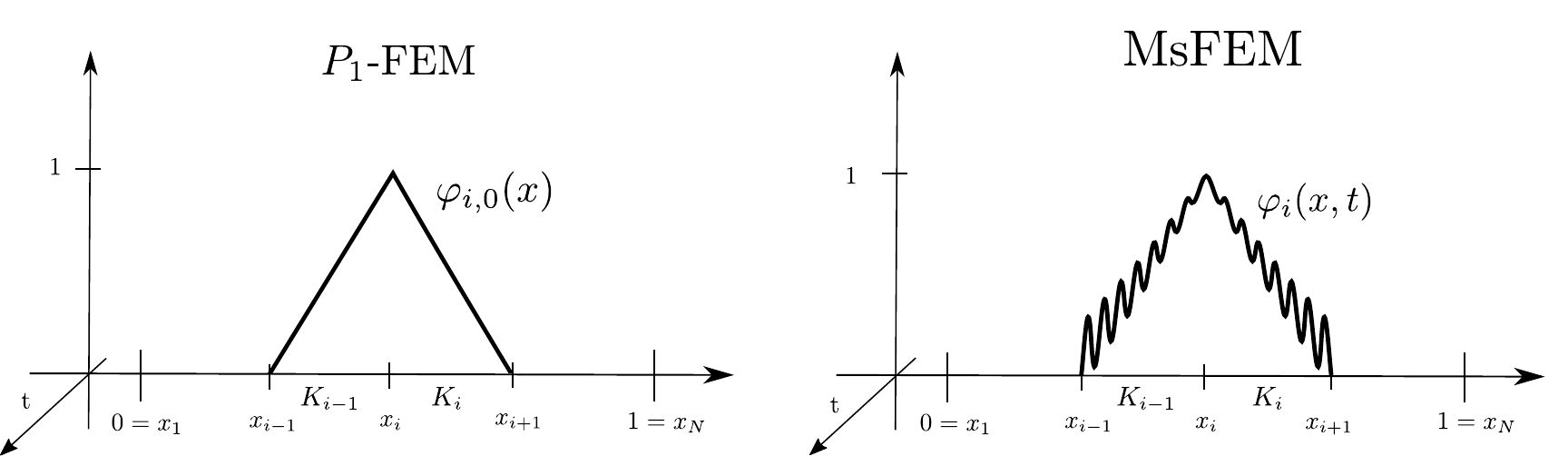}
  \caption{Illustration of basis functions. Standard FEM basis
    functions $\vphi_{i,0}(x,t)=\vphi_{i,0}(x)$ (left) have constant
    shape in time and serve as initial conditions for multiscale FEM
    basis functions $\vphi_i(x,t)$ (right) which does depend on
    time. The latter is glued together from parts that evolve in each
    adjacent cell according to the homogeneous equation. It carries the
    local structure of the flow in space and time.}
  \label{fig-basis}
\end{figure}

%% file: fig-characteristics.tex
\begin{figure}[b!]
  \centering
  \includegraphics[width=0.95\textwidth]{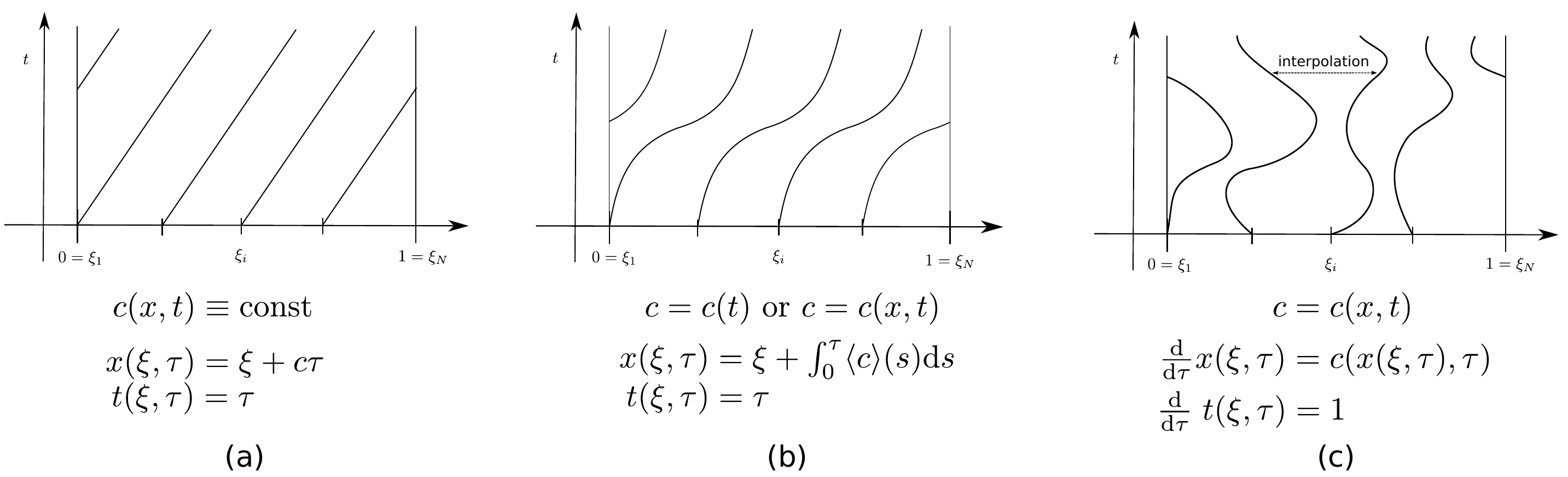}
  \caption{Illustration of different transforms. \textbf{(a)} If the
    background velocity $c(x,t)$ is constant transform~(\ref{eq:2-7a})
    gives a pure diffusion equation. \textbf{(b)} A spatially
    dependent $c$ is averaged over the entire domain and one follows
    the characteristics of the mean flow by
    transform~(\ref{eq:2-7b}). In the special case $c(x,t)=c(t)$ then
    the average equals $c(t)$ and the transformed equation is purely
    parabolic. \textbf{(c)} Transform~(\ref{eq:2-11}) follows the true
    characteristics that emerge from each coarse node $\xi_i$ and
    interpolates them linearly inside the coarse cells.}
  \label{fig-characteristics}
\end{figure}

%% file: sec-3.tex
\section{Numerical Tests}\label{s-3}

In this section we demonstrate our method on a number of examples. For
all examples presented here we use the same initial condition. In our
case this is a normalized Gaussian wave package, see
Figure~\ref{fig-initial_data} which is given by
\begin{equation}
  \label{eq:3-1}
  f(x) = \frac{1}{\sigma\sqrt{2\pi}}\exp \left( -\frac{(x-\nu)^2}{2\sigma^2}  \right)
\end{equation}
All data and all results will be set and shown from an Eulerian point
of view. Since no analytical reference solutions are available we will
compute a reference solution by a
high-resolution\input{s-3-fig-initial_data} standard FEM (750
elements, first order in space) in mean flow
coordinates~(\ref{eq:2-7b}), i.e., the advection dominance is already
reduced. For our multiscale FEM we will use $N=10$ coarse elements and
$N_f = 75$ fine elements in each coarse cell. Basis functions are
computed with a standard FEM. To show the advantage of our multiscale
FEM we will compare it to a low-resolution standard FEM that uses mean
flow coordinates like the reference solution and $N=10$ elements,
i.e., as many elements as the multiscale method uses on the coarse
scale. All simulations will be carried out on the time intervall
$[0,T]$ with $T=1$. We fix these parameters unless we explicitly say
something else.

\paragraph{Case 1.} We start with a simple setting. This test is a
direct generalization of existing
methods for elliptic equations. We
have time dependent coefficients given by
\begin{equation}
  \label{eq:3-2}
  \begin{split}
    c(t) & = 5\cos(10\pi t) \:, \\
    \mu(x,t) & = 5(t+1) \left( 0.01 + 0.0099\cos(2\pi k x)  \right) \:.
  \end{split}
\end{equation}
\input{s-3-fig-test-1-phase_space} Note that the background velocity
only depends on time. Therefore, characteristics can never cross and
hence transforms~(\ref{eq:2-7a}),~(\ref{eq:2-7b}) and~(\ref{eq:2-11})
are identical. Snapshots of the solution and the corresponding errors
to the reference solution are shown at four different time stamps in
Figure~\ref{s-3-fig-test-1-snapshots} and error graphs are shown in
Figure~\ref{s-3-fig-test-1-errors}. The characteristics induced by the
background velocity are shown in
Figure~\ref{s-3-fig-test-1-phase_space} and Table~\ref{table-test-1}
summarizes the errors in $L^2(I)$ and $L^\infty(I)$ for a number of
different $k$. Since the multiscale solution resolves the local
structure of the solution within coarse elements we also show the
error of the derivative, i.e., the $H^1(I)$-error. This can also
indirectly be seen in Figure~\ref{s-3-fig-test-1-snapshots} since the
error of the multiscale solution is smoother than the one of the
standard solution. The time step was taken to be $\delta
t=10^{-3}$. 

\vspace{1cm}
\input{s-3-table-test-1}

\paragraph{Case 2.} This test is to compare the two MsFEM versions
using either transform~(\ref{eq:2-7b}), i.e., the mean flow transform
(MF-MsFEM), or transform~(\ref{eq:2-11}) which we will refer to as the
characteristic MsFEM (Char-MsFEM). We test the influence of
oscillations in the background velocity that are either in the
resolved part of the spectrum or in the unresolved part. By this we
mean that we compare scenarios in which oscillations can or can not be
resolved by the coarse grid. The coefficients are given by
\begin{equation}
  \label{eq:3-3}
  \begin{split}
    c(x,t) & = 10 + \cos(2k\pi x) \:, \\
    \mu(x,t) & = 5(t+1)\left( 0.01 + 0.0099\cos(60\pi x)  \right) \:.
  \end{split}
\end{equation}
Here the advective term is on average stronger than the diffusive term
by a factor of approximately $5*10^2$ (comparing the mean values) and
by a factor of $10^5$ in regions where he diffusivity has a
minimum. Oscillations in the background velocity with $k<5$ can be
resolved on a grid with only $N=10$ elements whereas oscillations of
the diffusivity can not be
resolved. Figure~\ref{s-3-fig-test-2-snapshots} shows snapshots of the
solution for the case $k=3$ and the case $k=60$. Errors are shown in
Figure~\ref{s-3-fig-test-2-errors}. In the first case the Char-MsFEM
shows superior performance over both the standard method and the
MF-MsFEM. In the case of a larger frequency separation between the
mean velocity and the oscillatory part one can see that the MF-MsFEM
performs slightly better than the Char-MsFEM. Errors are summarized in
Table~\ref{table-test-2}. Note that the interesting scenario for us is
the one where $k<5$, i.e., the background velocity can be resolved by
the coarse grid since typically in climate simulations velocities are
taken from a coarse mesh and the diffusion parameter includes
parametrizations of processes that are not resolved. In such a
scenario our Char-MsFEM outperforms the MF-MsFEM whereas if $k>5$ we
observe a similar to slightly better performance of the MF-MsFEM.

\vspace{1cm}
\input{s-3-table-test-2}

\vspace{0.5cm}
\input{s-3-table-test-3} 

\vspace{0.5cm}
\paragraph{Case 3.} Here we demonstrate the advantage of the
Char-MsFEM in regimes dominated by advection where we have frequencies
in the background velocity that are resolved and that are unresolved
by the coarse grid. We do not show a comparison here to the MF-MsFEM
since this method fails if resolved frequencies are
apparent. Furthermore, we use an oscillatory diffusion coefficient as
we did in previous other tests. Coefficients for this test are given
by
\begin{equation}
  \label{eq:3-4}
  \begin{split}
    c(x,t) & = v + 1.5\cos(2\pi x) + 0.5\cos(60\pi x) \:, \\
    \mu(x,t) & = 5(t+1)\left( 0.01 + 0.0099\cos(50\pi x)  \right) \:.
  \end{split}
\end{equation}
Snapshots of the solution and errors for $v = 4$ are shown in
Figure~\ref{s-3-fig-test-3-snapshots} as well as a summary of errors
in Table~\ref{table-test-3}. From the snapshots it is clearly visible
that the standard FEM does not capture the correct large scale
behavior, i.e., it is too fast.

\input{s-3-fig-test-1-snapshots}
\input{s-3-fig-test-2-snapshots}
\input{s-3-fig-test-2-errors}

\input{s-3-fig-test-3-snapshots}

\paragraph{Case 4.} This test is to demonstrate consistency of the
Char-MsFEM in the limit of a large scale separation. This means that
we will demonstrate that our method yields a solution that correctly
follows the reference solution as the diffusion coefficient becomes
more and more oscillatory. It is not our aim to find the limiting
equation, i.e., a homogenized equation, and to demonstrate that our
numerical solution enjoys the correct behavior with respect to the
solution of the homogenized equation. Finding the correct homogenized
equation is often a difficult task, even in simple settings such as
purely diffusive ones. Additionally to the oscillating diffusive term
we will have a background velocity that is resolved by the coarse grid
since our motivation comes from climate simulations. There the
velocity data is usually taken from a coarse grid while the fine scale
part is hidden in the parametrizations. The coefficients that we use
for this case are
\begin{equation}
  \label{eq:3-5}
  \begin{split}
    c(x,t) & = (2t+0.5)  (3 + \cos(2\pi x) + \cos(60\pi x)) \:, \\
    \mu(x,t) & =  0.01 + 0.0099\cos(2k\pi x) \:.
  \end{split}
\end{equation}
Note that the standard solution also runs faster than the
Char-MsFEM. This effect is much more obvious than the wrong
representation of the fine scale diffusion on the coarse grid, see
Figure~\ref{s-3-fig-test-5-snapshots}. Similarly to case 1 but less
obvious the standard solution does not capture the oscillations in the
diffusive term well. This effect is less severe since the velocity
coefficient ``wipes'' these oscillations
out. Table~\ref{table-test-5-max} shows the relative deviation of the
maximum of the coarse standard solution and of the Char-MsFEM with
respect to the reference solution. Errors in $L^2$ and $L^\infty$ are
summarized in Table~\ref{table-test-5}.

\vspace{1cm}
\input{s-3-table-test-5}

\input{s-3-fig-test-5-snapshots}

\paragraph{Case 5.} We use a simple example to demonstrate that the
Char-MsFEM also works with non-vanishing external forcing. The
coefficients we choose are given by
\begin{equation}
  \label{eq:3-6}
  \begin{split}
    g(x) = & 0.015\sin(8\pi) \\
    c(x,t) & = (2t+0.5)  (3 + \cos(2\pi x) + \cos(60\pi x)) \:, \\
    \mu(x,t) & =  0.01 + 0.0099\cos(2k\pi x) \:.
  \end{split}
\end{equation}
Snapshots are shown in Figure~\ref{s-3-fig-test-4-snapshots}.
Remarkably, the Char-MsFEM also yields good results if we choose the
above forcing to act on a scale much finer than resolved by the coarse
grid. The resulting multiscale solution then of course does not
capture the fine scale forcing since the basis functions do not
(although they capture the fine scale diffusion) but we observed in
tests that it still runs correctly with the reference solution in
contrast to the solution of the standard FEM. However, we do not
strive to investigate external forcing on subgrid scales and therefore
pass on showing this.

\input{s-3-fig-test-4-snapshots}

\paragraph{A note on convergence.} In Case 4 we demonstrated that the
Char-MsFEM correctly captures fine scale features of the solution as
the scale separation in the diffusion coefficient grows. Here we show
indications for convergence of the Char-MsFEM as the coarse resolution
$H\rightarrow 0$. The classical multiscale finite element method for
elliptic problems as described in~\cite{Hou1997, Hou1999} is shown to
behave like a standard FEM as $H\rightarrow 0$. Here we see a similar
behavior but the Char-MsFEM convergences slowly. The data for the
tests is given by
\begin{equation}
  \label{eq:3-7}
  \begin{split}
    c(x,t) & = (2t+0.5)  (1.5 + \frac{\alpha}{2}\cos(2\pi x)) \:, \\
    \mu(x,t) & =  0.01 + 0.0099\cos(2k\pi x) \:.
  \end{split}
\end{equation}

The right-hand side $g$ is taken to be zero. Snapshots for the cases
$\alpha=0$ and $\alpha=1$ for $k=10$ are shown for different coarse
resolutions in Figure~\ref{s-3-fig-test-6-snapshots} as well as the
respective mean square and maximum errors in
Figure~\ref{s-3-fig-test-6-errors}. The influence of the smoothness of
the diffusion coefficient on the error taken with respect to a
high-resolution reference solution is summarized in
Table~\ref{table-test-6} for spatially constant background velocity as
well as for spatially varying velocity. The tests indicate that the
error is smaller the smoother the diffusion coefficient is relative to
the coarse grid, i.e., the method converges faster for smooth
diffusion. Clearly visible from Figure~\ref{s-3-fig-test-6-snapshots}
is that once the coefficients are resolved by the coarse grid the
Char-MsFEM solution is slower than the reference solution but this
effect vanishes as $H$ decreases. The time step was taken as
$\delta t=1/1000$ and the number of fine elements in each coarse block
was taken to be constant for all different $H$. However, note that
this regime is practically not relevant since here the standard FEM is
sufficient if the data is resolved by the coarse grid. The advantage
of the Char-MsFEM is at scales when data is not resolved by the coarse
grid.

%% file: s-3-fig-initial_data.tex
\begin{wrapfigure}[]{r}{0.5\textwidth}
  \begin{center}
    \vspace{-0.7cm}
      \includegraphics[width=0.48\textwidth]{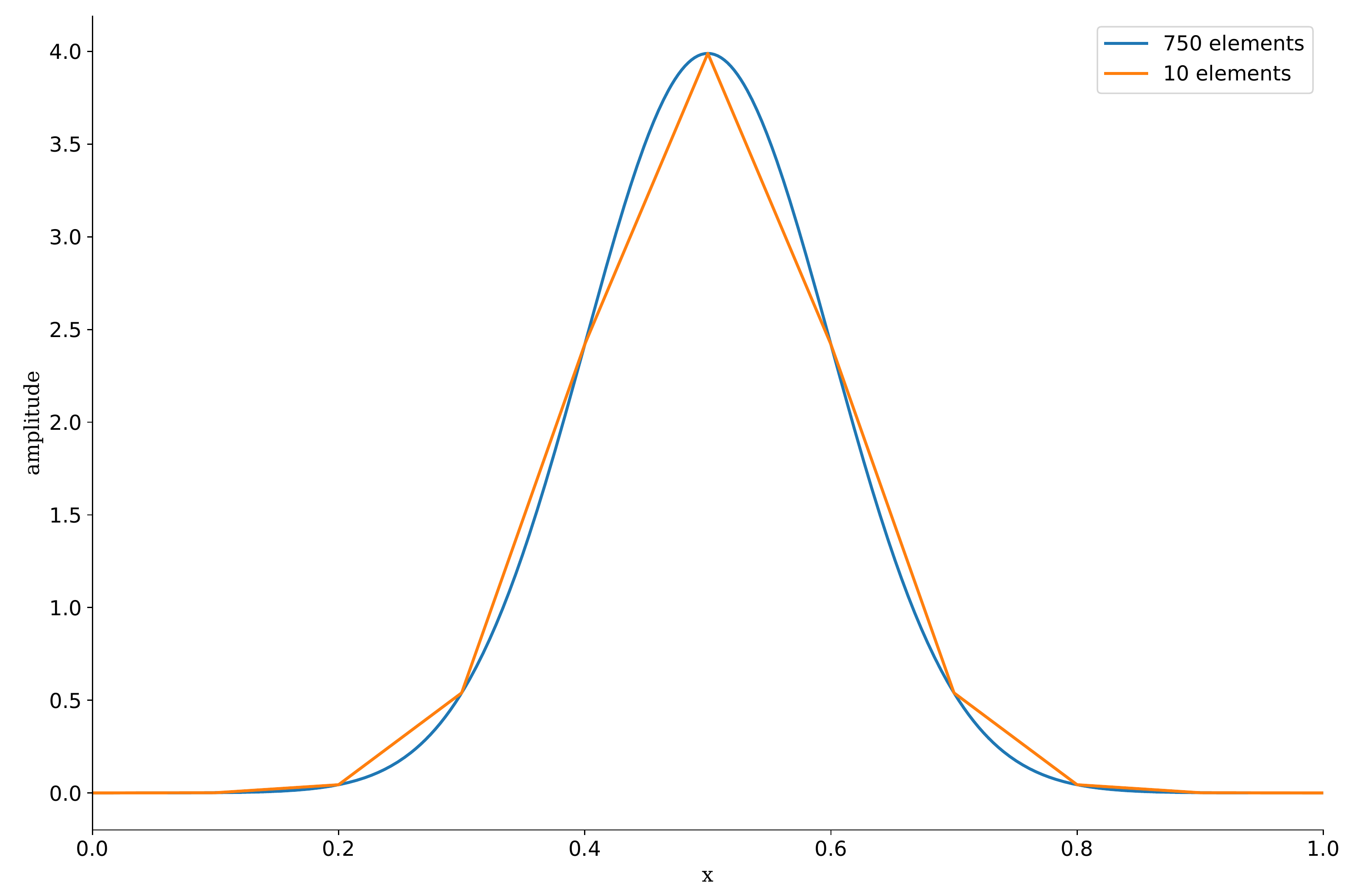}
  \end{center}
  \caption{Initial condition for all tests is given
    by~(\ref{eq:3-1}).}
  \label{fig-initial_data}
\end{wrapfigure}

%% file: s-3-fig-test-1-phase_space.tex
\begin{wrapfigure}[]{r}{0.4\textwidth}
  \begin{center}
    \vspace{-0.5cm}
      \includegraphics[width=0.38\textwidth]{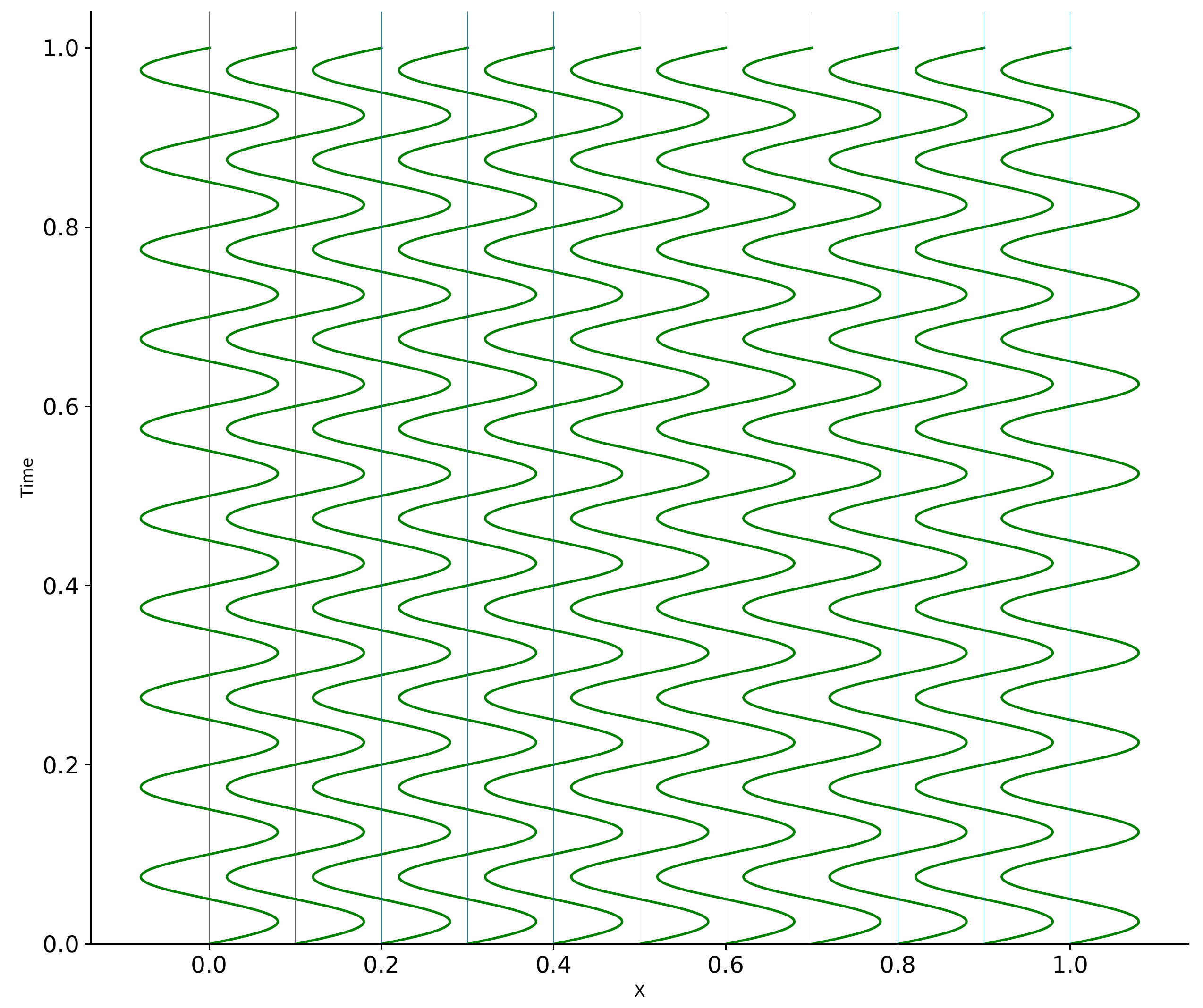}
  \end{center}
  \vspace{-0.5cm}
  \caption{Characteristics emerging at each coarse node induced by the
    velocity coefficient in~(\ref{eq:3-2}) with Characteristics
    originating at $\xi=0$ and $\xi=1$ are identical due to periodic
    boundary conditions.}
  \label{s-3-fig-test-1-phase_space}
\end{wrapfigure}

%% file: s-3-table-test-1.tex
\begin{table}[h]
  \begin{center}
    \begin{tabular}{ p{0.05\linewidth} | 
      p{0.075\textwidth} | p{0.075\textwidth} || 
      p{0.075\textwidth} | p{0.075\textwidth} ||
      p{0.075\textwidth} | p{0.075\textwidth} }

      & \multicolumn{2}{c||}{$L^2(I)$} & \multicolumn{2}{c||}{$L^\infty(I)$} & \multicolumn{2}{c|}{$H^1(I)$} \\

      & \multicolumn{1}{c|}{FEM} & \multicolumn{1}{c||}{MsFEM}
      & \multicolumn{1}{c|}{FEM} & \multicolumn{1}{c||}{MsFEM}
      & \multicolumn{1}{c|}{FEM} & \multicolumn{1}{c|}{MsFEM} \\ \hline

      \multicolumn{1}{c|}{$k=15$} 
      & \multicolumn{1}{c|}{$1.036\cdot10^{-4}$} & \multicolumn{1}{c||}{$1.848\cdot10^{-5}$} 
      & \multicolumn{1}{c|}{$1.141\cdot10^{-4}$} & \multicolumn{1}{c||}{$2.258\cdot10^{-5}$} 
      & \multicolumn{1}{c|}{$5.622\cdot10^{-4}$} & \multicolumn{1}{c|}{$2.160\cdot10^{-4}$}  \\ \hline

      \multicolumn{1}{c|}{$k=30$} 
      & \multicolumn{1}{c|}{$1.770\cdot 10^{-4}$} & \multicolumn{1}{c||}{$1.437\cdot10^{-5}$} 
      & \multicolumn{1}{c|}{$1.809\cdot 10^{-4}$} & \multicolumn{1}{c||}{$1.785\cdot 10^{-5}$} 
      & \multicolumn{1}{c|}{$6.329\cdot 10^{-4}$} & \multicolumn{1}{c|}{$1.706\cdot 10^{-4}$}  \\ \hline
      
      \multicolumn{1}{c|}{$k=45$} 
      & \multicolumn{1}{c|}{$2.054\cdot10^{-4}$} & \multicolumn{1}{c||}{$1.160\cdot10^{-5}$} 
      & \multicolumn{1}{c|}{$2.031\cdot10^{-4}$} & \multicolumn{1}{c||}{$1.460\cdot10^{-5}$} 
      & \multicolumn{1}{c|}{$6.299\cdot10^{-4}$} & \multicolumn{1}{c|}{$1.556\cdot10^{-4}$}  \\ \hline

      \multicolumn{1}{c|}{$k=60$} 
      & \multicolumn{1}{c|}{$2.370\cdot10^{-4}$} & \multicolumn{1}{c||}{$1.143\cdot10^{-5}$} 
      & \multicolumn{1}{c|}{$2.293\cdot10^{-4}$} & \multicolumn{1}{c||}{$1.405\cdot10^{-5}$} 
      & \multicolumn{1}{c|}{$6.523\cdot10^{-4}$} & \multicolumn{1}{c|}{$1.887\cdot10^{-4}$}  \\ \hline

      \hline
    \end{tabular}
  \end{center}
  \caption{Relative error of case 1 at $T=1$. Coefficients are given by~(\ref{eq:3-2}).}
  \label{table-test-1}
\end{table}

%% file: s-3-table-test-2.tex
\begin{table}[H]
  \begin{center}
    \begin{tabular}{ p{0.05\linewidth} | 
      p{0.075\textwidth} | p{0.075\textwidth} | p{0.075\textwidth} }

      \multicolumn{1}{c}{} & \multicolumn{3}{c}{$L^2(I)$} \\
      & \multicolumn{1}{c|}{FEM} & \multicolumn{1}{c|}{MF-MsFEM} & \multicolumn{1}{c|}{Char-MsFEM} \\ \hline
      \multicolumn{1}{c|}{$k=3$} 
      & \multicolumn{1}{c|}{$4.052\cdot10^{-5}$} & \multicolumn{1}{c|}{$4.284\cdot10^{-5}$} & \multicolumn{1}{c|}{$1.038\cdot10^{-5}$}  \\ \hline
      \multicolumn{1}{c|}{$k=15$}  
      & \multicolumn{1}{c|}{$4.080\cdot 10^{-5}$} & \multicolumn{1}{c|}{$3.205\cdot10^{-6}$} & \multicolumn{1}{c|}{$6.787\cdot 10^{-6}$}  \\ \hline
      \multicolumn{1}{c|}{$k=30$} 
      & \multicolumn{1}{c|}{$1.171\cdot10^{-4}$} & \multicolumn{1}{c|}{$1.634\cdot10^{-5}$} & \multicolumn{1}{c|}{$7.543\cdot10^{-6}$}  \\ \hline
      \multicolumn{1}{c|}{$k=60$} 
      & \multicolumn{1}{c|}{$4.765\cdot10^{-5}$} & \multicolumn{1}{c|}{$7.046\cdot10^{-6}$} & \multicolumn{1}{c|}{$2.868\cdot10^{-5}$} \\ \hline

      \multicolumn{1}{c}{} & \multicolumn{3}{c}{$L^\infty(I)$} \\
      & \multicolumn{1}{c|}{FEM} & \multicolumn{1}{c|}{MF-MsFEM} & \multicolumn{1}{c|}{Char-MsFEM} \\ \hline
      \multicolumn{1}{c|}{$k=3$} 
      & \multicolumn{1}{c|}{$5.258\cdot10^{-5}$} & \multicolumn{1}{c|}{$5.404\cdot10^{-5}$} & \multicolumn{1}{c|}{$1.743\cdot10^{-5}$}  \\ \hline
      \multicolumn{1}{c|}{$k=15$} 
      & \multicolumn{1}{c|}{$5.071\cdot10^{-5}$} & \multicolumn{1}{c|}{$6.818\cdot10^{-6}$} & \multicolumn{1}{c|}{$1.542\cdot10^{-5}$}  \\ \hline
      \multicolumn{1}{c|}{$k=30$} 
      & \multicolumn{1}{c|}{$1.385\cdot10^{-4}$} & \multicolumn{1}{c|}{$2.332\cdot10^{-5}$} & \multicolumn{1}{c|}{$1.087\cdot10^{-5}$}  \\ \hline
      \multicolumn{1}{c|}{$k=60$} 
      & \multicolumn{1}{c|}{$5.969\cdot10^{-5}$} & \multicolumn{1}{c|}{$9.408\cdot10^{-6}$} & \multicolumn{1}{c|}{$3.863\cdot10^{-5}$} \\ \hline

      \multicolumn{1}{c}{}  & \multicolumn{3}{c}{$H^1(I)$} \\
      & \multicolumn{1}{c|}{FEM} & \multicolumn{1}{c|}{MF-MsFEM} & \multicolumn{1}{c|}{Char-MsFEM} \\ \hline
      \multicolumn{1}{c|}{$k=3$} 
      & \multicolumn{1}{c|}{$5.201\cdot10^{-4}$} & \multicolumn{1}{c|}{$4.108\cdot10^{-4}$} & \multicolumn{1}{c|}{$2.987\cdot10^{-4}$}  \\ \hline
      \multicolumn{1}{c|}{$k=15$} 
      & \multicolumn{1}{c|}{$5.185\cdot10^{-4}$} & \multicolumn{1}{c|}{$2.539\cdot10^{-4}$} & \multicolumn{1}{c|}{$2.609\cdot10^{-4}$}  \\ \hline
      \multicolumn{1}{c|}{$k=30$} 
      & \multicolumn{1}{c|}{$7.323\cdot10^{-4}$} & \multicolumn{1}{c|}{$3.467\cdot10^{-4}$} & \multicolumn{1}{c|}{$2.665\cdot10^{-4}$}  \\ \hline
      \multicolumn{1}{c|}{$k=60$} 
      & \multicolumn{1}{c|}{$5.558\cdot10^{-4}$} & \multicolumn{1}{c|}{$2.459\cdot10^{-4}$} & \multicolumn{1}{c|}{$5.210\cdot10^{-4}$}  \\ \hline

      \hline
    \end{tabular}
  \end{center}
  \caption{Relative error of case 2 at $T=1$. Coefficients are given by~(\ref{eq:3-3}).}
  \label{table-test-2}
\end{table}

%% file: s-3-table-test-3.tex
\begin{table}[h!]
  \begin{center}
    \begin{tabular}{ p{0.05\linewidth} | 
      p{0.075\textwidth} | p{0.075\textwidth} || 
      p{0.075\textwidth} | p{0.075\textwidth} }

      & \multicolumn{2}{c||}{$L^2(I)$} & \multicolumn{2}{c|}{$L^\infty(I)$} \\

      & \multicolumn{1}{c|}{FEM} & \multicolumn{1}{c||}{Char-MsFEM}
      & \multicolumn{1}{c|}{FEM} & \multicolumn{1}{c|}{Char-MsFEM} \\ \hline

      \multicolumn{1}{c|}{$v=4$} 
      & \multicolumn{1}{c|}{$2.688\cdot10^{-4}$} & \multicolumn{1}{c||}{$5.635\cdot10^{-5}$} 
      & \multicolumn{1}{c|}{$2.692\cdot10^{-4}$} & \multicolumn{1}{c|}{$6.779\cdot10^{-5}$}  \\ \hline

      \multicolumn{1}{c|}{$v=8$} 
      & \multicolumn{1}{c|}{$1.722\cdot 10^{-4}$} & \multicolumn{1}{c||}{$4.731\cdot10^{-5}$} 
      & \multicolumn{1}{c|}{$1.684\cdot 10^{-4}$} & \multicolumn{1}{c|}{$6.077\cdot 10^{-5}$}  \\ \hline
      
      \multicolumn{1}{c|}{$v=12$} 
      & \multicolumn{1}{c|}{$1.286\cdot10^{-4}$} & \multicolumn{1}{c||}{$2.213\cdot10^{-5}$} 
      & \multicolumn{1}{c|}{$1.281\cdot10^{-4}$} & \multicolumn{1}{c|}{$2.848\cdot10^{-5}$}  \\ \hline

      \multicolumn{1}{c|}{$v=16$} 
      & \multicolumn{1}{c|}{$9.850\cdot10^{-4}$} & \multicolumn{1}{c||}{$2.736\cdot10^{-5}$} 
      & \multicolumn{1}{c|}{$9.924\cdot10^{-4}$} & \multicolumn{1}{c|}{$3.874\cdot10^{-5}$}  \\ \hline

      \hline
    \end{tabular}
  \end{center}
  \caption{Relative error of case 3 at $T=1$. Coefficients are given by~(\ref{eq:3-4}).}
  \label{table-test-3}
\end{table}

%% file: s-3-fig-test-1-snapshots.tex
\clearpage
\thispagestyle{empty}

\begin{figure*}[t!]
  \centering
  \includegraphics[width=1\textwidth]{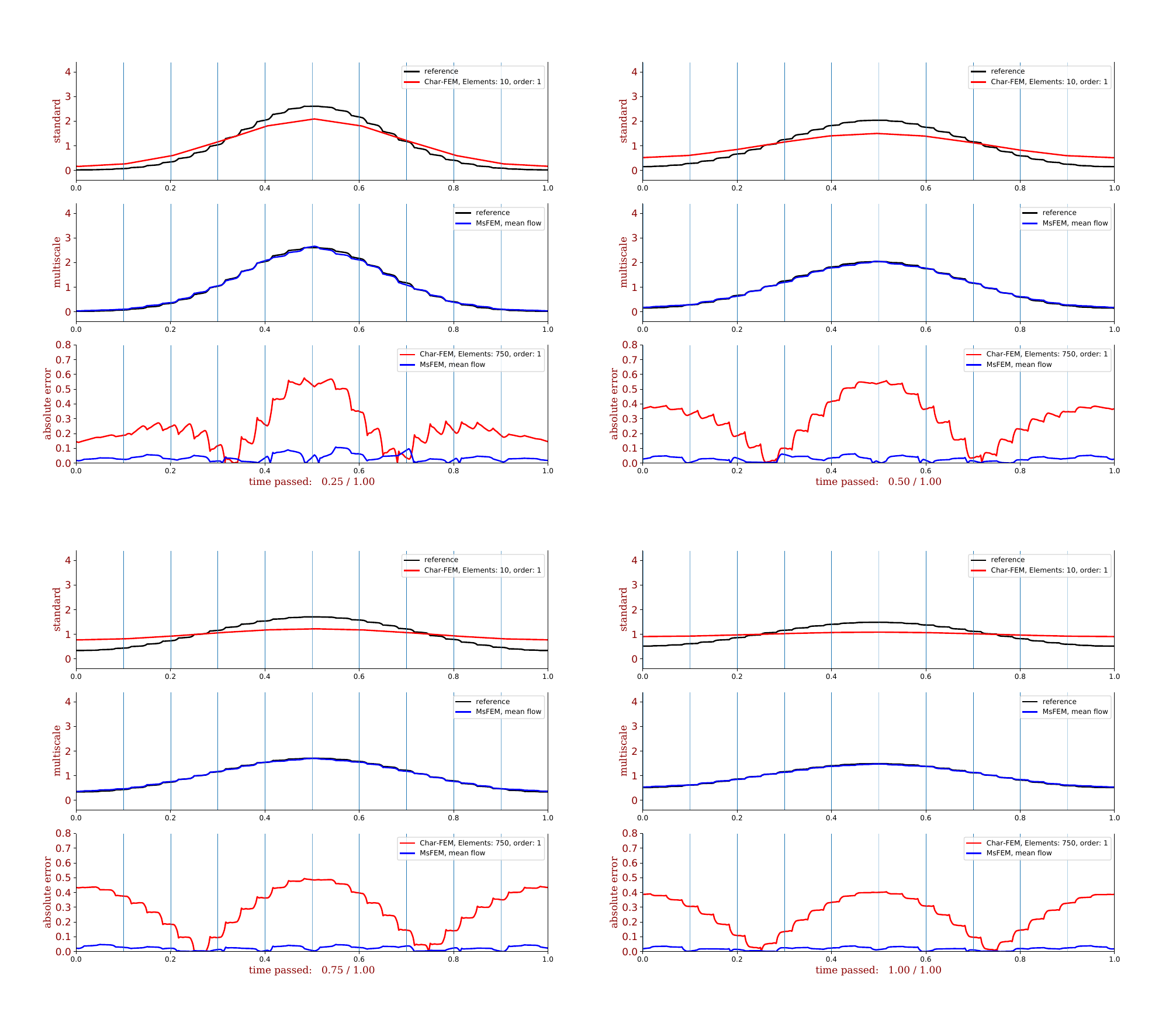}
  \caption[]{Snapshots of the solution of case 1 with data given
    by~(\ref{eq:3-2}) at time stamps $t=0.25$, $t=0.5$, $=0.75$ and
    $T=1$ for $k=30$. The upper graphs show a low-resolution standard
    FEM (red) compared to a high-resolution standard FEM reference
    solution (black). The second row of each snapshot shows our MsFEM
    (blue) compared to the reference solution. The third row of each
    snapshot shows the absolute difference of the low-resolution FEM
    solution (red) and of the MsFEM solution (blue) to the reference
    solution. The absolute error of the MsFEM is clearly lower and
    smoother that the error of the standard FEM.}
  \label{s-3-fig-test-1-snapshots}
\end{figure*}

\begin{figure*}[h!]
  \centering
  \includegraphics[width=0.99\textwidth]{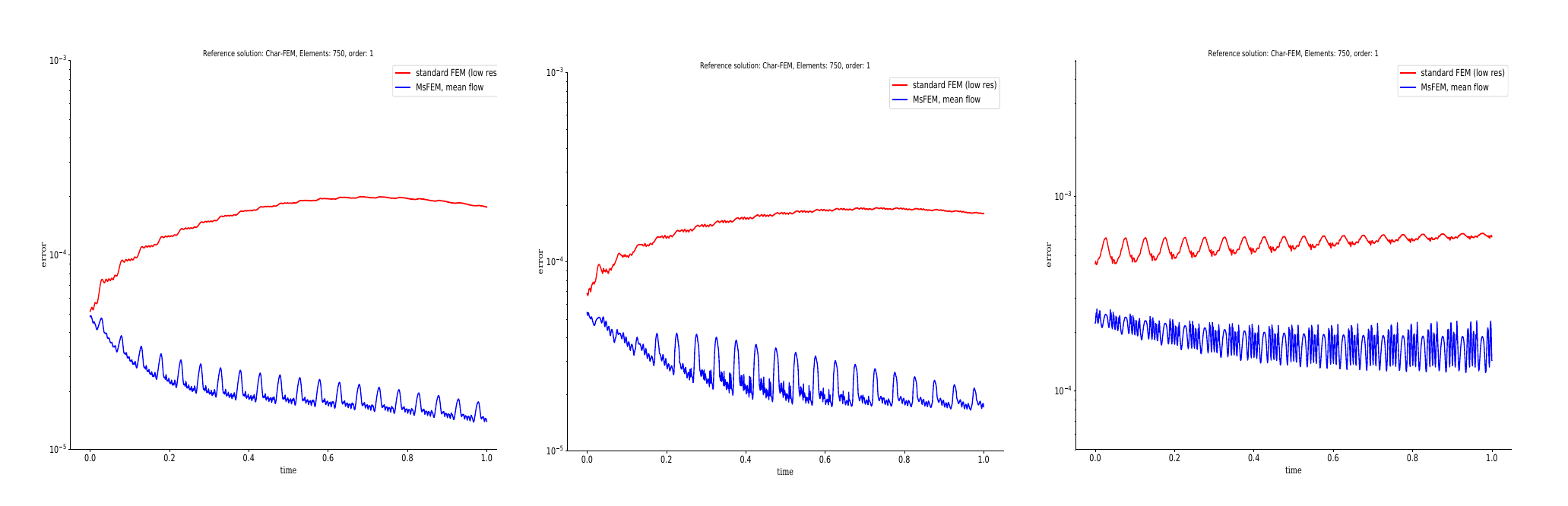}
  \caption[]{Errors of the solution of case 1 with data given
    by~(\ref{eq:3-2}) for $k=30$. Errors of a low-resolution standard
    FEM (red) are compared to errors of our MsFEM (blue) with respect
    to a high-resolution reference solution. The $L^2$-error and the
    $L^\infty$-error of the MsFEM are approximately one order of
    magnitude better than the errors of the standard FEM. The error of
    the derivative of the MsFEM is approximately half an order of
    magnitude better than the error of the standard FEM.}
  \label{s-3-fig-test-1-errors}
\end{figure*}

\clearpage

%% file: s-3-fig-test-2-snapshots.tex
\clearpage
\thispagestyle{empty}

\begin{figure*}[t!]
  \centering
\vspace{-2cm}
  \includegraphics[width=0.95\textwidth]{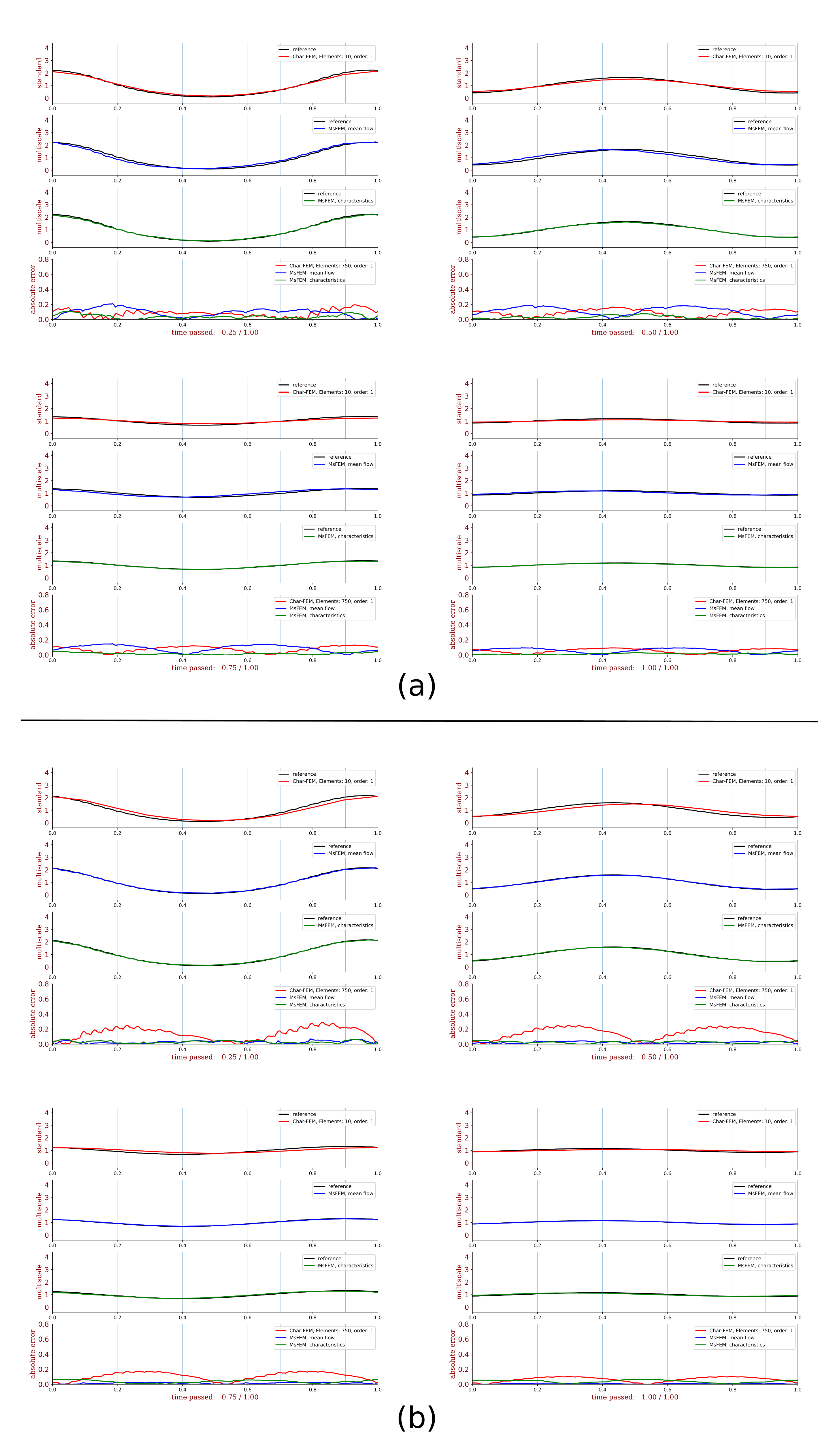}
  \caption[]{Comparison of snapshots of the solution of case 2 with
    data given by~(\ref{eq:3-3}) at time stamps $t=0.25$, $t=0.5$,
    $=0.75$ and $T=1$ for \textbf{(a)} $k=3$ and \textbf{(b)}
    $k=60$. In each image the first row represents a comparison of
    either a standard FEM (red), the MF-MsFEM (blue) or a Char-MsFEM
    (green) to a reference solution (black). The respective absolute
    errors are shown in the last row of each image.}
  \label{s-3-fig-test-2-snapshots}
\end{figure*}

\clearpage

%% file: s-3-fig-test-2-errors.tex
\clearpage
\thispagestyle{empty}

\begin{figure*}[t!]
  \centering
  \includegraphics[width=0.95\textwidth]{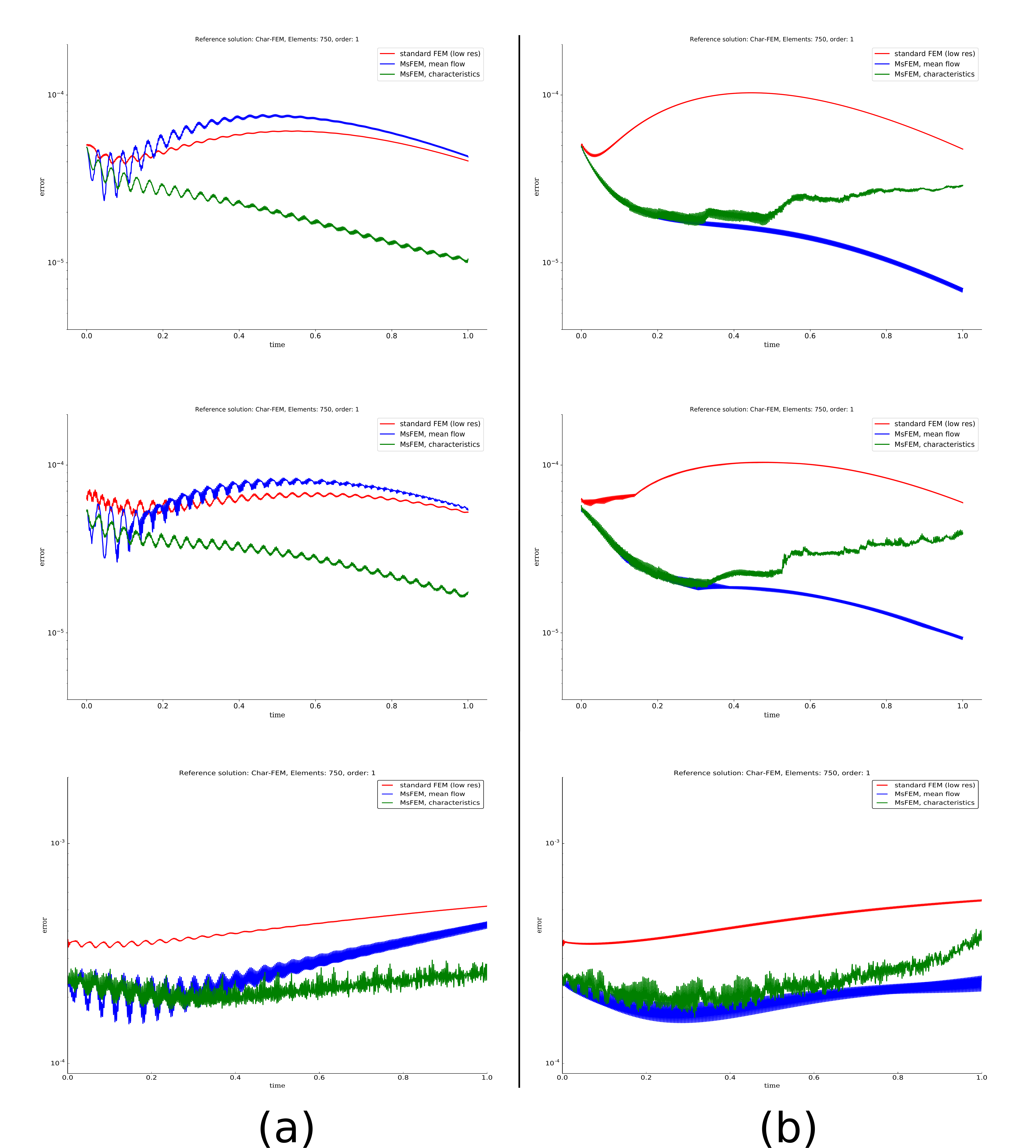}
  \caption[]{Errors of the solution of case 2 with data given
    by~(\ref{eq:3-3}) for \textbf{(a)} $k=3$ and \textbf{(b)} $k=60$.}
  \label{s-3-fig-test-2-errors}
\end{figure*}

\clearpage

%% file: s-3-fig-test-3-snapshots.tex
\clearpage
\thispagestyle{empty}

\begin{figure*}[t!]
  \centering
  \vspace{-3cm}
  \includegraphics[width=1\textwidth]{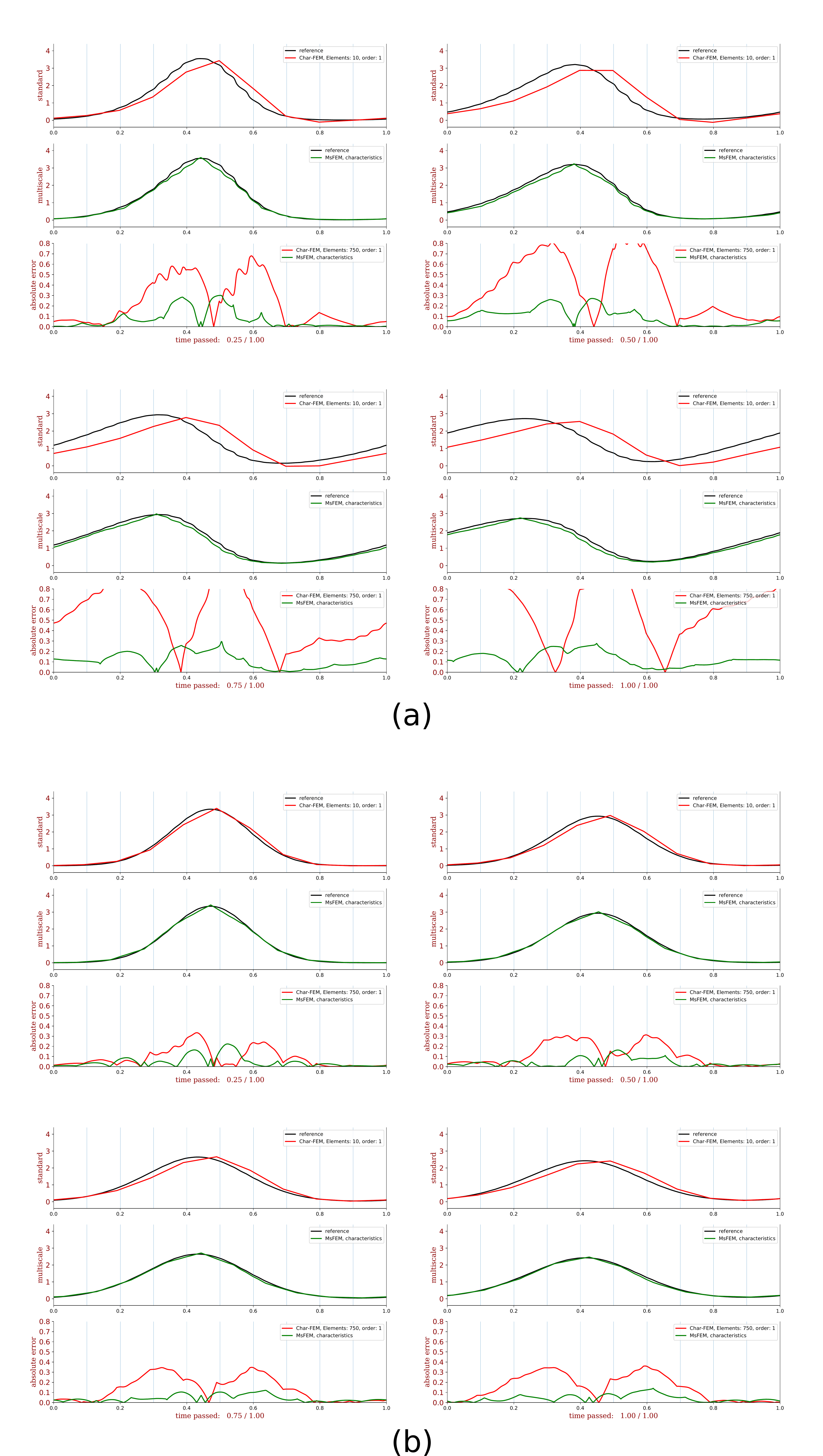}
  \caption[]{Comparison of snapshots of the solution of case 3 with
    data given by~(\ref{eq:3-4}) at time stamps $t=0.25$, $t=0.5$,
    $=0.75$ and $T=1$ for \textbf{(a)} $v=4$ and \textbf{(b)}
    $v=16$. Each image contains a comparison of a standard FEM (red)
    and the Char-MsFEM (green) to a reference solution (black). The
    absolute errors with respective colors are shown in the last row
    of each image.}
  \label{s-3-fig-test-3-snapshots}
\end{figure*}

\clearpage

%% file: s-3-table-test-5.tex
\begin{table}[h!]
  \begin{center}
    \begin{tabular}{ p{0.05\linewidth} | 
      p{0.075\textwidth} | p{0.075\textwidth} }

      & \multicolumn{1}{c|}{FEM} & \multicolumn{1}{c|}{Char-MsFEM} \\ \hline

      \multicolumn{1}{c|}{$k=15$} 
      & \multicolumn{1}{c|}{$3.66\cdot10^{-2}$} & \multicolumn{1}{c|}{$8.94\cdot10^{-3}$}  \\ \hline

      \multicolumn{1}{c|}{$k=30$} 
      & \multicolumn{1}{c|}{$5.67\cdot10^{-2}$} & \multicolumn{1}{c|}{$9.63\cdot10^{-3}$}  \\ \hline

      \multicolumn{1}{c|}{$k=60$} 
      & \multicolumn{1}{c|}{$7.80\cdot10^{-2}$} & \multicolumn{1}{c|}{$2.43\cdot10^{-2}$}  \\ \hline

      \multicolumn{1}{c|}{$k=90$} 
      & \multicolumn{1}{c|}{$8.52\cdot10^{-2}$} & \multicolumn{1}{c|}{$3.08\cdot10^{-2}$}  \\ \hline

      \multicolumn{1}{c|}{$k=140$} 
      & \multicolumn{1}{c|}{$1.09\cdot10^{-1}$} & \multicolumn{1}{c|}{$1.50\cdot10^{-2}$}  \\ \hline

      \multicolumn{1}{c|}{$k=200$} 
      & \multicolumn{1}{c|}{$1.47\cdot10^{-1}$} & \multicolumn{1}{c|}{$1.74\cdot10^{-2}$}  \\ \hline
      \hline
    \end{tabular}
  \end{center}
  \caption{Relative deviation of the maximum of case 4 at $T=1$. Coefficients are given by~(\ref{eq:3-5}).}
  \label{table-test-5-max}
\end{table}

\vspace{0.5cm}

\begin{table}[h!]
  \begin{center}
    \begin{tabular}{ p{0.05\linewidth} | 
      p{0.075\textwidth} | p{0.075\textwidth} || 
      p{0.075\textwidth} | p{0.075\textwidth} }

      & \multicolumn{2}{c||}{$L^2(I)$} & \multicolumn{2}{c|}{$L^\infty(I)$} \\

      & \multicolumn{1}{c|}{FEM} & \multicolumn{1}{c||}{Char-MsFEM}
      & \multicolumn{1}{c|}{FEM} & \multicolumn{1}{c|}{Char-MsFEM} \\ \hline

      \multicolumn{1}{c|}{$k=15$} 
      & \multicolumn{1}{c|}{$2.137\cdot10^{-4}$} & \multicolumn{1}{c||}{$5.418\cdot10^{-5}$} 
      & \multicolumn{1}{c|}{$2.204\cdot10^{-4}$} & \multicolumn{1}{c|}{$5.168\cdot10^{-5}$}  \\ \hline

      \multicolumn{1}{c|}{$k=30$} 
      & \multicolumn{1}{c|}{$2.173\cdot10^{-4}$} & \multicolumn{1}{c||}{$2.233\cdot10^{-5}$} 
      & \multicolumn{1}{c|}{$2.196\cdot10^{-4}$} & \multicolumn{1}{c|}{$5.142\cdot10^{-5}$}  \\ \hline

      \multicolumn{1}{c|}{$k=60$} 
      & \multicolumn{1}{c|}{$2.220\cdot10^{-4}$} & \multicolumn{1}{c||}{$4.898\cdot10^{-5}$} 
      & \multicolumn{1}{c|}{$2.189\cdot10^{-4}$} & \multicolumn{1}{c|}{$4.731\cdot10^{-5}$}  \\ \hline

      \multicolumn{1}{c|}{$k=90$} 
      & \multicolumn{1}{c|}{$2.240\cdot10^{-4}$} & \multicolumn{1}{c||}{$4.585\cdot10^{-5}$} 
      & \multicolumn{1}{c|}{$2.203\cdot10^{-4}$} & \multicolumn{1}{c|}{$4.755\cdot10^{-5}$}  \\ \hline

      \multicolumn{1}{c|}{$k=140$} 
      & \multicolumn{1}{c|}{$2.308\cdot10^{-4}$} & \multicolumn{1}{c||}{$4.481\cdot10^{-5}$} 
      & \multicolumn{1}{c|}{$2.211\cdot10^{-4}$} & \multicolumn{1}{c|}{$5.028\cdot10^{-5}$}  \\ \hline

      \multicolumn{1}{c|}{$k=200$} 
      & \multicolumn{1}{c|}{$2.420\cdot10^{-4}$} & \multicolumn{1}{c||}{$4.933\cdot10^{-5}$} 
      & \multicolumn{1}{c|}{$2.215\cdot10^{-4}$} & \multicolumn{1}{c|}{$6.083\cdot10^{-5}$}  \\ \hline
      \hline
    \end{tabular}
  \end{center}
  \caption{Relative error of case 4 at $T=1$. Coefficients are given by~(\ref{eq:3-5}).}
  \label{table-test-5}
\end{table}

%% file: s-3-fig-test-5-snapshots.tex
\clearpage
\thispagestyle{empty}

\begin{figure*}[t!]
  \centering
  \vspace{-3cm}
  \includegraphics[width=0.86\textwidth]{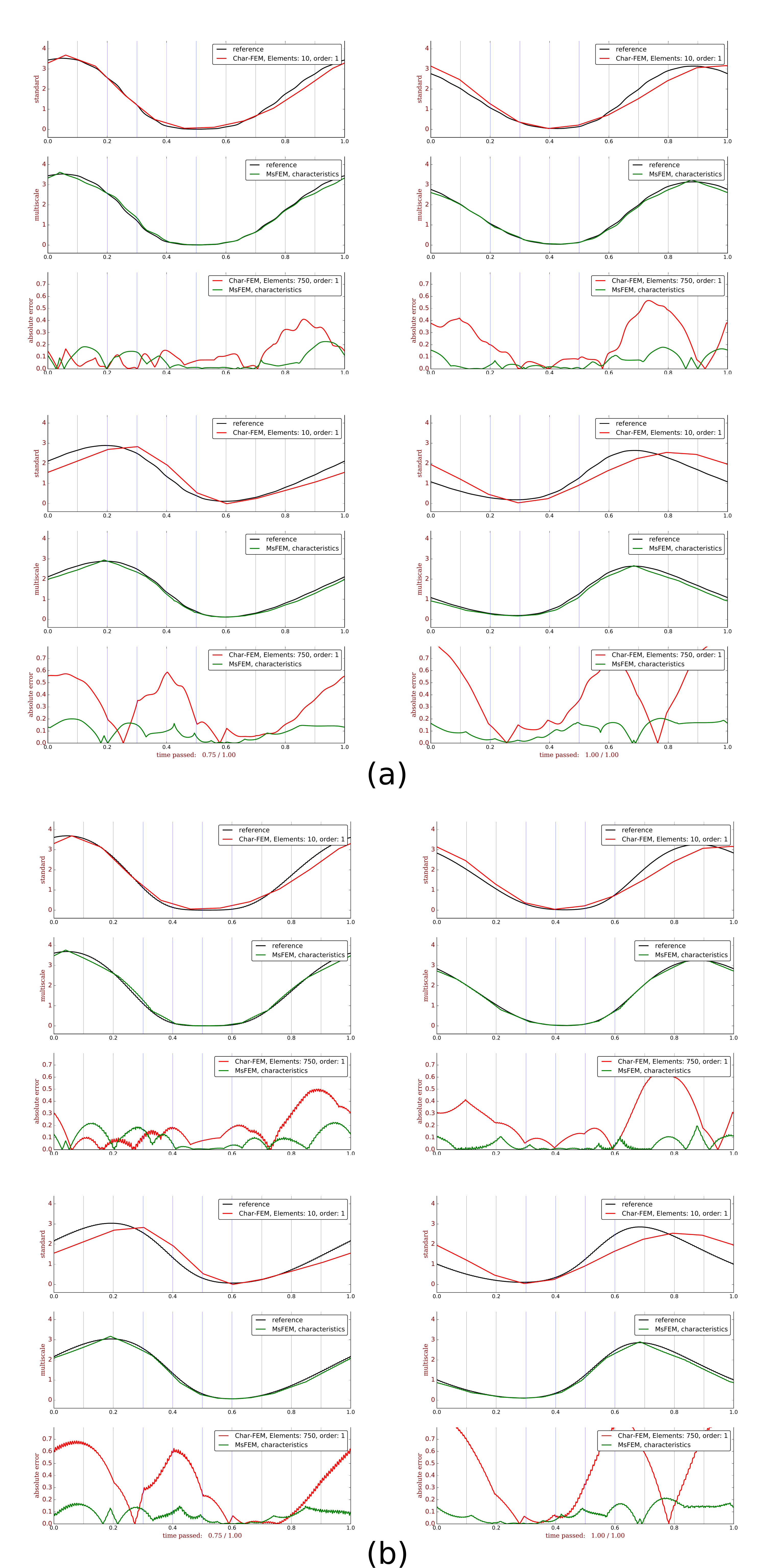}
  \caption[]{Comparison of snapshots of the solution of case 4 with
    data given by~(\ref{eq:3-5}) at time stamps $t=0.25$, $t=0.5$,
    $=0.75$ and $T=1$ for \textbf{(a)} $k=15$ and \textbf{(b)}
    $k=140$. Each image contains a comparison of a standard FEM (red)
    and the Char-MsFEM (green) to a reference solution (black). The
    absolute errors with respective colors are shown in the last row
    of each image.}
  \label{s-3-fig-test-5-snapshots}
\end{figure*}

\clearpage

%% file: s-3-fig-test-4-snapshots.tex
\begin{figure*}[t!]
  \centering
  \vspace{-3cm}
  \includegraphics[width=1\textwidth]{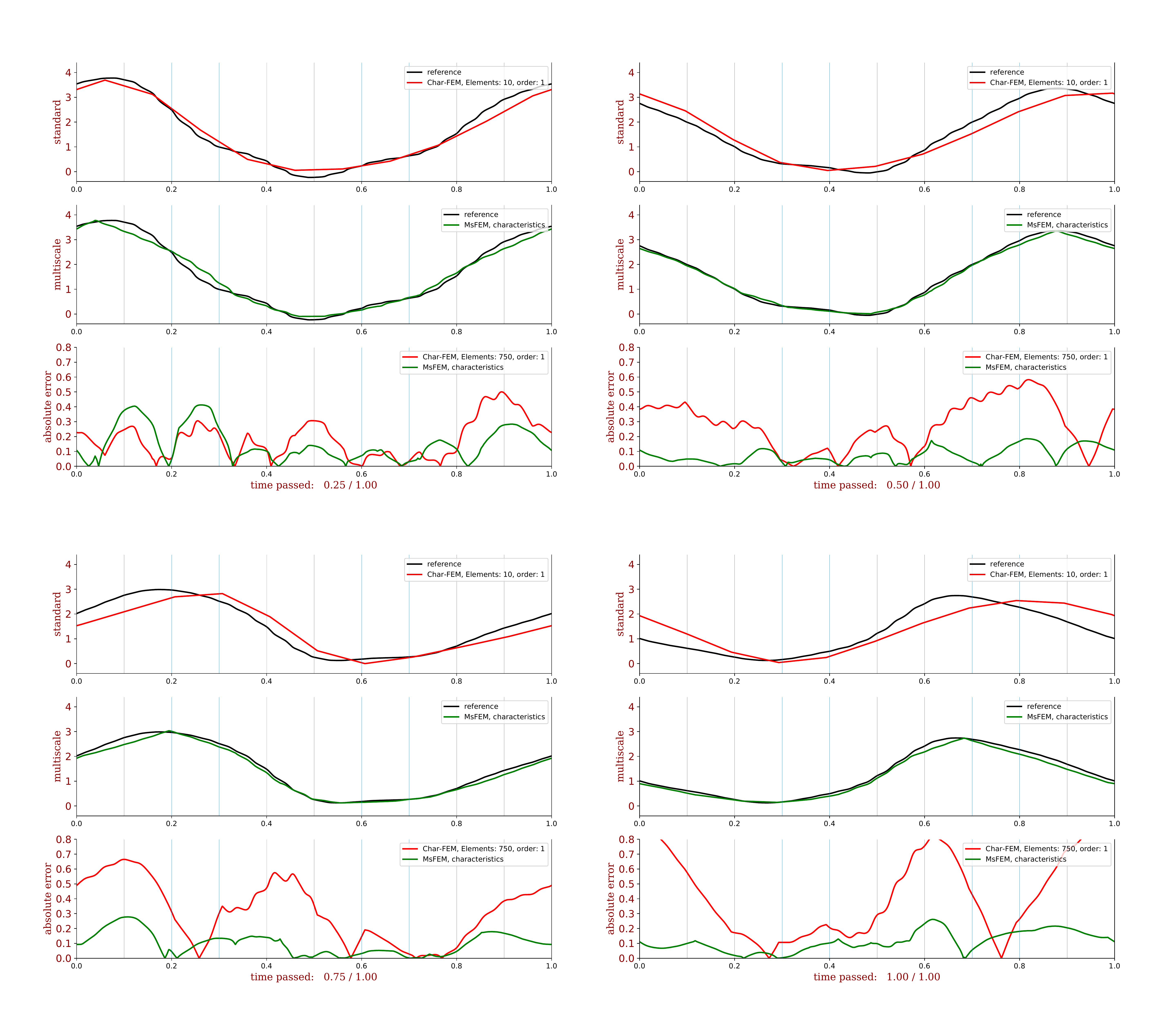}
  \caption[]{Comparison of snapshots of the solution of case 5 with
    data given by~(\ref{eq:3-6}) at time stamps $t=0.25$, $t=0.5$,
    $=0.75$ and $T=1$. Each image contains a comparison of a standard
    FEM (red) and the Char-MsFEM (green) to a reference solution
    (black). The absolute errors with respective colors are shown in
    the last row of each image.}
  \label{s-3-fig-test-4-snapshots}
\end{figure*}

%% file: sec-4.tex
\section{Discussion}\label{s-4}

We have suggested a flexible numerical multiscale method based on a
finite element formulation that uses a set of non-polynomial basis
functions to solve an advection-diffusion equation with a smooth,
well-behaved and dominant background velocity and a heterogeneous
diffusive term. The model problem is motivated by problems that arise
in long-term climate simulations where prognostic variables like the
background velocity live on a coarse grid while parametrized subgrid
processes are modelled in the diffusive part.

The numerical method is motivated by a classical MsFEM suggested by
the porous media community. We substantially modified the ideas to
account for scenarios relevant to climate simulations. This
modification is based on transforms suggested by the background
velocity. In particular these transforms aim to make the dominant
advective term milder. Our suggested multiscale methods are not overly
complicated and consist of an offline phase in which many small local
problems for the basis functions (with Dirichlet boundary conditions
on curves suggested by the transform) are solved and an online phase
that uses the basis in a finite element like framework in the global
formulation. The offline phase is suitable for massive parallelization
where as the online phase is fast. Also the global form is not
restricted to use a finite element framework. It can potentially be
used with finite volumes or discontinuous Galerkin formulations.

In tests we showed an improvement of numerical solutions using our new
technique compared to standard solutions on the same coarse
grid. Since no analytical solutions are available we compare our
solutions to high-resolution standard FEM solutions. Tests indicate
that our multiscale solutions enjoy the correct upscaling behavior in
the limit of fast oscillations in the diffusion, i.e., the faster the
diffusion oscillates relative to the coarse grid the more advantageous
are our multiscale methods over standard methods.

Our multiscale techniques show good results but unfortunately their
scope is limited to velocity fields that have a dominant mean value
for all times. This has to do with the dynamics of the different
characteristics. The less dominant the mean of the background velocity
the closer characteristics emerging at coarse grid boundaries can come
in short time. In particular time-independent zeros in the background
velocity represent stable and unstable attractors of the advective
part. Such features can make coarse blocks collapse and make the
method unstable.  

Furthermore, a direct generalization to higher dimensions is quite
difficult. This is due to richer dynamics in the characteristics even
if the velocity field is solenoidal. Also, our methods fail if the
advection-diffusion equation is in conservation form. The reason is
that a conservative transport term of the form $(cu)_x$ decomposes
into an advective part $cu_x$ and a reactive part $c_xu$ that
compensates for non-divergence of the background velocity $c$. Such
reaction terms we found in tests enjoy a different upscaling mechanism
and can not be treated with our method.

However, although our method does not work in all cases we showed that
it works in relevant scenarios and identified difficulties in their
generalization. Being aware of such difficulties can facilitate the
choice of modified ingredients to account for generalizations and for
the development

\clearpage
\thispagestyle{empty}
\input{s-3-fig-test-6-errors}
\input{s-3-table-test-6}
\clearpage

\input{s-3-fig-test-6-snapshots}

\noindent of suitable multiscale methods. For example, the reader may
think of a combination of semi-Lagrangian methods to deal with the
problem of converging or crossing characteristics and our multiscale
methods. This is left for future research. With this work we hope to
spark a debate on Galerkin based multiscale methods in the climate
simulation community and to contribute to their success.

%% file: s-3-fig-test-6-errors.tex
\begin{figure*}[h!]
  \centering
  \includegraphics[width=1\textwidth]{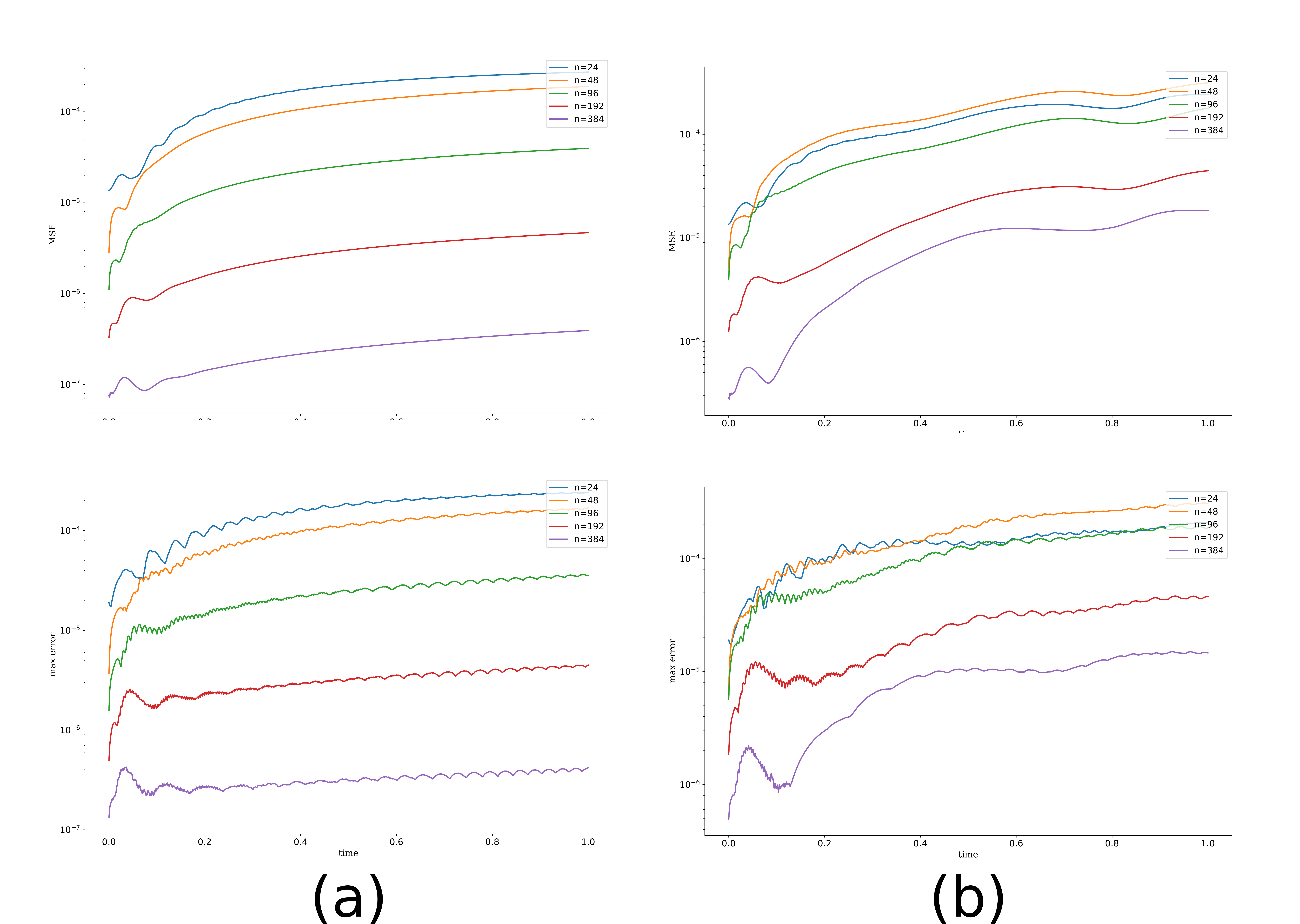}
  \caption[]{Relative errors of the solution with data given
    by~(\ref{eq:3-7}) for $k=10$ over time for \textbf{(a)} $\alpha=0$
    and \textbf{(b)} $\alpha=1$.}
  \label{s-3-fig-test-6-errors}
\end{figure*}

%% file: s-3-table-test-6.tex
\begin{table}[H]
\centering
\setlength\tabcolsep{4pt}
\begin{minipage}{0.48\textwidth}
  \centering
  \begin{tabular}{cccc}
    
    \multicolumn{1}{c|}{$k$} & \multicolumn{1}{c|}{$N$} & \multicolumn{1}{c|}{$\norm{e}{L^2}$} & \multicolumn{1}{c|}{$\norm{e}{L^\infty}$}  \\ \hline
    
    \multicolumn{1}{c|}{$10$} & \multicolumn{1}{c|}{$24$} & \multicolumn{1}{c|}{$2.732\cdot10^{-4}$} & \multicolumn{1}{c|}{$2.435\cdot10^{-4}$}  \\
    \multicolumn{1}{c|}{} & \multicolumn{1}{c|}{$48$} & \multicolumn{1}{c|}{$1.898\cdot10^{-4}$} & \multicolumn{1}{c|}{$1.692\cdot10^{-4}$}  \\
    \multicolumn{1}{c|}{} & \multicolumn{1}{c|}{$96$} & \multicolumn{1}{c|}{$3.963\cdot10^{-5}$} & \multicolumn{1}{c|}{$3.569\cdot10^{-5}$}  \\
    \multicolumn{1}{c|}{} & \multicolumn{1}{c|}{$192$} & \multicolumn{1}{c|}{$4.688\cdot10^{-6}$} & \multicolumn{1}{c|}{$4.474\cdot10^{-6}$}  \\
    \multicolumn{1}{c|}{} & \multicolumn{1}{c|}{$384$} & \multicolumn{1}{c|}{$3.939\cdot10^{-7}$} & \multicolumn{1}{c|}{$4.204\cdot10^{-7}$}  \\ \hline
    
    \multicolumn{1}{c|}{$15$} & \multicolumn{1}{c|}{$24$} & \multicolumn{1}{c|}{$9.603\cdot10^{-5}$} & \multicolumn{1}{c|}{$9.129\cdot10^{-5}$}  \\
    \multicolumn{1}{c|}{} & \multicolumn{1}{c|}{$48$} & \multicolumn{1}{c|}{$3.400\cdot10^{-4}$} & \multicolumn{1}{c|}{$2.911\cdot10^{-4}$}  \\
    \multicolumn{1}{c|}{} & \multicolumn{1}{c|}{$96$} & \multicolumn{1}{c|}{$1.141\cdot10^{-4}$} & \multicolumn{1}{c|}{$9.920\cdot10^{-5}$}  \\
    \multicolumn{1}{c|}{} & \multicolumn{1}{c|}{$192$} & \multicolumn{1}{c|}{$1.731\cdot10^{-5}$} & \multicolumn{1}{c|}{$1.566\cdot10^{-5}$}  \\
    \multicolumn{1}{c|}{} & \multicolumn{1}{c|}{$384$} & \multicolumn{1}{c|}{$1.455\cdot10^{-6}$} & \multicolumn{1}{c|}{$1.368\cdot10^{-6}$}  \\ \hline
    
    \multicolumn{1}{c|}{$20$} & \multicolumn{1}{c|}{$24$} & \multicolumn{1}{c|}{$2.629\cdot10^{-5}$} & \multicolumn{1}{c|}{$2.730\cdot10^{-5}$}  \\
    \multicolumn{1}{c|}{} & \multicolumn{1}{c|}{$48$} & \multicolumn{1}{c|}{$3.819\cdot10^{-4}$} & \multicolumn{1}{c|}{$3.250\cdot10^{-4}$}  \\
    \multicolumn{1}{c|}{} & \multicolumn{1}{c|}{$96$} & \multicolumn{1}{c|}{$1.986\cdot10^{-4}$} & \multicolumn{1}{c|}{$1.686\cdot10^{-4}$}  \\
    \multicolumn{1}{c|}{} & \multicolumn{1}{c|}{$192$} & \multicolumn{1}{c|}{$4.110\cdot10^{-5}$} & \multicolumn{1}{c|}{$3.616\cdot10^{-5}$}  \\
    \multicolumn{1}{c|}{} & \multicolumn{1}{c|}{$384$} & \multicolumn{1}{c|}{$4.104\cdot10^{-6}$} & \multicolumn{1}{c|}{$3.702\cdot10^{-6}$}  \\ \hline
    
    \hline
  \end{tabular}
  \subcaption{$\alpha = 0$}
\end{minipage}%
\hfill
\begin{minipage}{0.48\textwidth}
  \centering
  \begin{tabular}{cccc}
    
    \multicolumn{1}{c|}{$k$} & \multicolumn{1}{c|}{$N$} & \multicolumn{1}{c|}{$\norm{e}{L^2}$} & \multicolumn{1}{c|}{$\norm{e}{L^\infty}$}  \\ \hline
    
    \multicolumn{1}{c|}{$10$} & \multicolumn{1}{c|}{$24$} & \multicolumn{1}{c|}{$2.409\cdot10^{-4}$} & \multicolumn{1}{c|}{$2.105\cdot10^{-4}$}  \\
    \multicolumn{1}{c|}{} & \multicolumn{1}{c|}{$48$} & \multicolumn{1}{c|}{$3.140\cdot10^{-4}$} & \multicolumn{1}{c|}{$3.096\cdot10^{-4}$}  \\
    \multicolumn{1}{c|}{} & \multicolumn{1}{c|}{$96$} & \multicolumn{1}{c|}{$1.783\cdot10^{-4}$} & \multicolumn{1}{c|}{$1.888\cdot10^{-4}$}  \\
    \multicolumn{1}{c|}{} & \multicolumn{1}{c|}{$192$} & \multicolumn{1}{c|}{$4.447\cdot10^{-5}$} & \multicolumn{1}{c|}{$4.640\cdot10^{-5}$}  \\
    \multicolumn{1}{c|}{} & \multicolumn{1}{c|}{$384$} & \multicolumn{1}{c|}{$1.828\cdot10^{-5}$} & \multicolumn{1}{c|}{$1.465\cdot10^{-5}$}  \\ \hline
    
    \multicolumn{1}{c|}{$15$} & \multicolumn{1}{c|}{$24$} & \multicolumn{1}{c|}{$1.647\cdot10^{-4}$} & \multicolumn{1}{c|}{$1.335\cdot10^{-4}$}  \\
    \multicolumn{1}{c|}{} & \multicolumn{1}{c|}{$48$} & \multicolumn{1}{c|}{$5.162\cdot10^{-4}$} & \multicolumn{1}{c|}{$4.487\cdot10^{-4}$}  \\
    \multicolumn{1}{c|}{} & \multicolumn{1}{c|}{$96$} & \multicolumn{1}{c|}{$4.159\cdot10^{-4}$} & \multicolumn{1}{c|}{$4.128\cdot10^{-4}$}  \\
    \multicolumn{1}{c|}{} & \multicolumn{1}{c|}{$192$} & \multicolumn{1}{c|}{$9.032\cdot10^{-5}$} & \multicolumn{1}{c|}{$9.610\cdot10^{-5}$}  \\
    \multicolumn{1}{c|}{} & \multicolumn{1}{c|}{$384$} & \multicolumn{1}{c|}{$2.001\cdot10^{-5}$} & \multicolumn{1}{c|}{$1.965\cdot10^{-5}$}  \\ \hline
    
    \multicolumn{1}{c|}{$20$} & \multicolumn{1}{c|}{$24$} & \multicolumn{1}{c|}{$7.963\cdot10^{-5}$} & \multicolumn{1}{c|}{$6.090\cdot10^{-5}$}  \\
    \multicolumn{1}{c|}{} & \multicolumn{1}{c|}{$48$} & \multicolumn{1}{c|}{$6.185\cdot10^{-4}$} & \multicolumn{1}{c|}{$5.368\cdot10^{-4}$}  \\
    \multicolumn{1}{c|}{} & \multicolumn{1}{c|}{$96$} & \multicolumn{1}{c|}{$6.595\cdot10^{-4}$} & \multicolumn{1}{c|}{$6.122\cdot10^{-4}$}  \\
    \multicolumn{1}{c|}{} & \multicolumn{1}{c|}{$192$} & \multicolumn{1}{c|}{$1.817\cdot10^{-4}$} & \multicolumn{1}{c|}{$1.847\cdot10^{-4}$}  \\
    \multicolumn{1}{c|}{} & \multicolumn{1}{c|}{$384$} & \multicolumn{1}{c|}{$2.930\cdot10^{-5}$} & \multicolumn{1}{c|}{$3.160\cdot10^{-5}$}  \\ \hline
    
    \hline
  \end{tabular}
  \subcaption{$\alpha = 1$}
\end{minipage}%
\caption{Relative errors of the numerical solution for data given by~(\ref{eq:3-7}) for 
  \textbf{(a)} $\alpha=0$ and \textbf{(b)} $\alpha = 1$ as $H$ becomes small, 
  i.e., $N$ becomes large. The table shows the influence of increasing frequency
  of the diffusion coefficient.}
\label{table-test-6}
\end{table}

%% file: s-3-fig-test-6-snapshots.tex
\clearpage
\thispagestyle{empty}

\begin{figure*}[t!]
  \centering
  \vspace{-1cm}
  \includegraphics[width=1\textwidth]{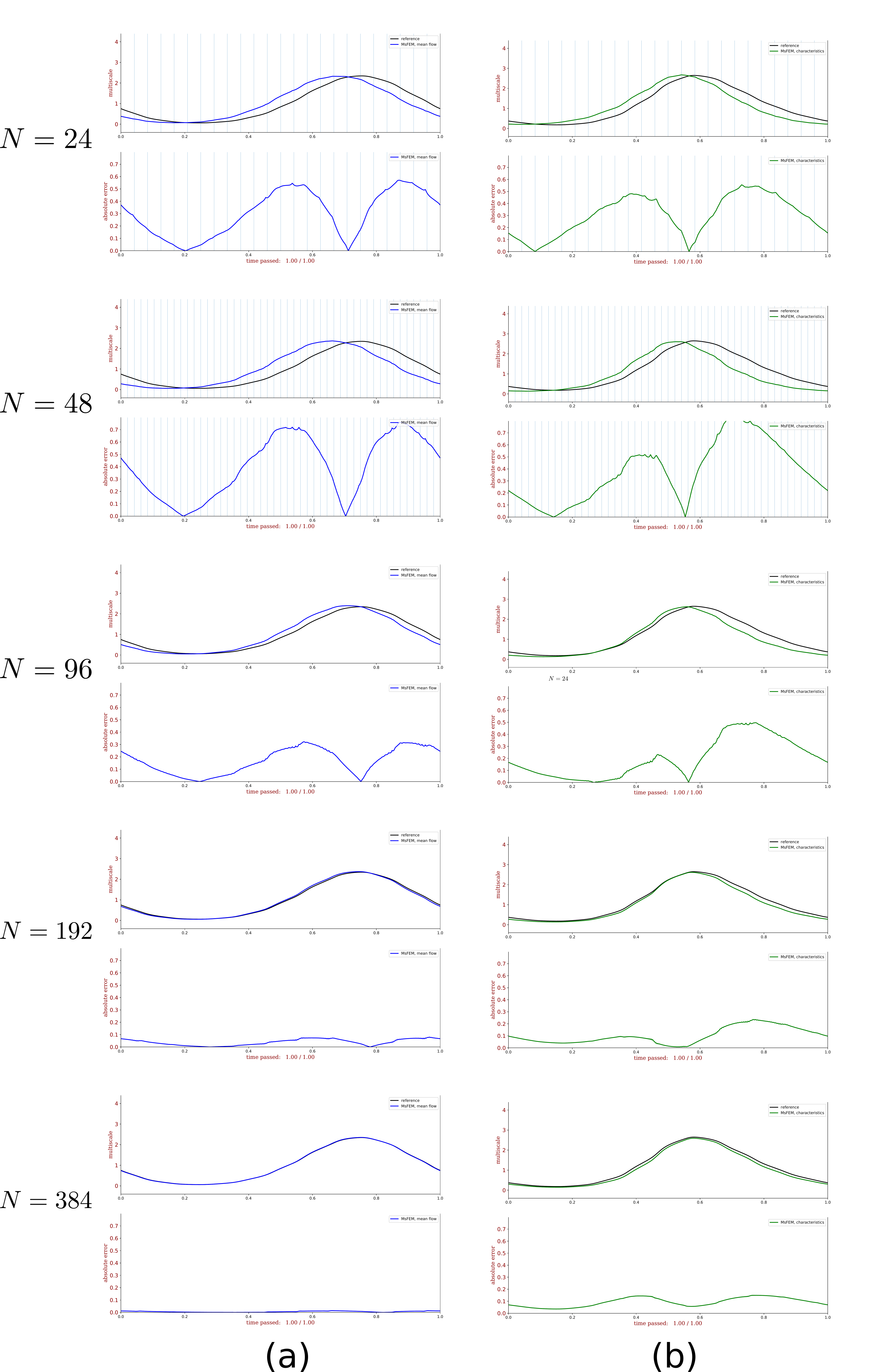}
  \caption[]{Comparison of snapshots of the solution with data given
    by~(\ref{eq:3-7}) at time $T=1$ and $k=10$ for \textbf{(a)}
    $\alpha=0$ and \textbf{(b)} $\alpha=1$. The multiscale solution is
    compared to a reference solution (black). The absolute errors are
    shown in the second row of each image.}
  \label{s-3-fig-test-6-snapshots}
\end{figure*}

\clearpage